\documentclass[a4paper,10pt]{article}

\title{PENALISATION OF THE SYMMETRIC RANDOM WALK\\by several functions of the supremum}
\author{DEBS PIERRE}
\date{}
\usepackage{amsthm}
\usepackage{a4wide,amsfonts,amsmath,latexsym,amssymb,euscript,graphicx,dsfont,units,color}
\usepackage[frenchb,english]{babel}
 
\newcommand{\p}{\mathds{P}}
\newcommand{\e}{\mathds{E}}

\theoremstyle{plain}     
\newtheorem{thm}{Theorem}[section]     
     
\newtheorem{lem}[thm]{Lemma}  
          
\newtheorem{rem}[thm]{Remark} 
     
\theoremstyle{definition}     
    
\begin{document}

\maketitle

\begin{abstract}
Call $(\Omega,\mathcal F_\infty,\p,X,\mathcal F)$ the canonical space
for the standard random walk on $\mathbb Z$. Thus,
$\Omega$ denotes the set of paths  $\phi:\mathbb N
\to\mathbb Z$ such that  ${|\phi(n+1)-\phi(n)|=1}$, 
$X=\left( X_n,n\geq0\right)$ is the canonical coordinate process 
on $\Omega$; $\mathcal F =\left( \mathcal F _n,n\geq0\right) $ 
is the natural filtration of $X$, $\mathcal F _\infty$
the $\sigma$-field $\bigvee_{n\geq0}\mathcal F _n$,
and $\p_0$
the probability on $\left( \Omega,\mathcal F _\infty\right)$ 
such that under $\p_0$, $X$ is the
standard random walk started from $0$, i.e., 
$\p_0\left( X_{n+1}=j\,\vert\,X_n=i\right)=\frac{1}{2}$ when $|j-i|=1$. 

Let $G:\mathbb N \times \Omega \rightarrow {\mathbb R}^+$ 
be a positive, adapted functional. For several types of functionals
$G$, we show the existence of a positive $\mathcal F$-martingale
$(M_n,\ n\geq0)$ such that, for all $n$  and all
$\Lambda_n\in\mathcal F _n$, 
$$\frac{\e_0[\mathds{1}_{\Lambda_n}G_p]}{\e_0[G_p]}
\quad\longrightarrow\quad\e_0[\mathds{1}_{\Lambda_n}M_n]\qquad
\hbox{when \ $p\rightarrow\infty$}\;. $$
Thus, there exists a probability $Q$ on 
$(\Omega,\mathcal F _\infty)$ such that
$Q(\Lambda_n)=\e_0[\mathds{1}_{\Lambda_n}M_n]$
for all $\Lambda_n\in\mathcal F _n$. We describe the behavior of
the process $(\Omega,X,\mathcal F )$ under~$Q$.\\

We study here four kinds of $G$:\\
.$G_p$ is a function of $S_p$ where $S_p$ is the unilateral supremum of $X$.\\
.$G_p$ is a function of $S_{g_p}$ where $g_p$ is the last 0 at the left of $p$.\\
.$G_p$ is a function of $S_{d_p}$ where $d_p$ is the first 0 at the right of $p$.\\
.$G_p$ is a function of $S^*_{g_p}$ where $S^*_p$ is the bilateral supremum of $X$.\\
.$G_p$ is a function of $S^*_{p}$ .\\
A similar study has been realized for other kinds of $G$ (cf \cite{D:1}).

\end{abstract}
\section{Introduction}
Let $\bigl\lbrace \Omega, {( X_t, \mathcal F_t)}_{t\geq 0} ,
\mathcal F_\infty,\p_x  \bigr\rbrace$ be the canonical
one-dimensional Brownian motion. For several types of positive functionals 
$\Gamma:\mathbb R^+ \times \Omega \rightarrow {\mathbb R}^+$,
B. Roynette, P. Vallois and M. Yor show in \cite{RVY:1} 
that, for fixed  $s$ and for all $\Lambda_s\in\mathcal F _s$,
$$\lim_{t\to\infty}\frac{\e_x[\mathds{1}_{\Lambda_s}
\Gamma_t]}{\e_x[\Gamma_t]}$$
exists and has the form $\e_x[\mathds{1}_{\Lambda_s}M_s^x]$, 
where $(M_s^x,s\geq0)$ is a positive martingale. 
This enables them to define a probability $Q_x$ on
$(\Omega,\mathcal F _\infty)$ by:
$$\forall\Lambda_s\in\mathcal F _s\qquad Q_x(\Lambda_s)
=\e_x[\mathds{1}_{\Lambda_s}M_s^x]\;;$$
moreover, they precisely describe the behavior of 
the canonical process $X$
under $Q_x$. They do this for numerous functionals
$\Gamma$, for instance a function of the one-sided maximum,
or of the local time, or of the age of the current excursion
(cf.~\cite{RVY:1}, \cite{RVY:7}).\\
We have already studied a discrete analogue of their results in \cite{D:1}.
More precisely, let $\Omega$ denote the set of all functions 
$\phi$ from $\mathbb N$ to $\mathbb Z$ such that 
${|\phi(n+1)-\phi(n)|=1}$, let $X=\left( X_n,n\geq0\right)$ be the 
process of coordinates on that space, 
$\mathcal F=(\mathcal F _n,n\geq0)$ the canonical filtration, 
$\mathcal F _\infty$ the $\sigma$-field 
$\bigvee_{n\geq0}\mathcal F _n$, and $\p_x$ $(x\in\mathbb N)$ the
family of probabilities on $\left( \Omega,\mathcal F _\infty\right)$ 
such that under $\p_x$, $X$ is the standard random walk started at~$x$. 
For notational simplicity, we often write $\p$ for $\p_0$.
Our aim is to establish that for several types
of positive, adapted functionals
$G:\mathbb N\times \Omega\rightarrow\mathbb N$,\\
i) for each $n\geq0$ and each $\Lambda_n\in\mathcal F _n$,
\begin{equation*}
\frac{\e_0[\mathds{1}_{\Lambda_n}G_p]}{\e_0[G_p]},
\end{equation*}
tends to a limit when $p$ tends to infinity;\\
ii) this limit is equal to
$\e_0[\mathds{1}_{\Lambda_n}M_n]$, for some
$\mathcal F$-martingale $M$ such that $M_0=1$.\\

Call $Q({\Lambda_n})$ this limit. Like the continuous case, $Q$ describes a probability on $\left(\Omega,\mathcal F_\infty\right)$ by :
$$\forall n\geq0,\forall \Lambda_n\in\mathcal F_n,\,Q(\Lambda_n)=\e_0[\mathds{1}_{\Lambda_n}M_n],$$ 
and we also study the process $X$ under $Q$.\\
A better definition of the principle of penalisation, for instance proof of existence and unicity,  can be found in the introduction and the first part of \cite{D:1}.\\

In this paper, $G$ essentially depends on two functions $\varphi :\mathbb N\rightarrow\mathbb R^+$ and $\phi :\mathbb N\rightarrow\mathbb R^+$ such that :
\begin{equation}\label{conditions}
\sum_{k\geq0}\varphi\left( k\right)=1,\,\phi\left( x\right):=\sum_{k=x}^{\infty}\varphi\left( k\right).
\end{equation}
The following result comes from \cite{D:1}, and is not proved in the following paper. Here,  $G$ is a function of the one-sided maximum, i.e. 
$G_p=\varphi(S_p)$, 
where $S_p:=\sup\left\lbrace X_k,k\leq p\right\rbrace$. 
We establish :
\begin{thm}\label{thmmax}
\begin{enumerate}
\item
\begin{enumerate}
\item For each $n\geq0$ and each 
$\Lambda_{n}\in{\mathcal{F}}_n$, one has
\begin{equation*}
\lim_{p\rightarrow\infty}\frac{\e[\mathds{1}_{\Lambda_n}
\varphi(S_p)]}{\e[\varphi(S_p)]}=\e[\mathds{1}_{\Lambda_n}
M_n^{\varphi}]\;,
\end{equation*}
where $M_n^{\varphi}:=\varphi(S_n)(S_n-X_n)+\phi(S_n)$.\\
\item $(M_n^{\varphi},n\geq0)$ is a positive martingale, with
$M^\varphi_0=1$, non 
uniformly integrable; in fact, $M_n^\varphi$ tends a.s. to $0$ when
$n\rightarrow\infty$.\\
\end{enumerate}
\item Call $Q^\varphi$ the probability  on 
$\left( \Omega,\mathcal F_\infty\right)$ characterized by
\begin{equation*}
\forall n\in\mathbb N, \Lambda_n\in\mathcal F _n,\quad 
Q^\varphi(\Lambda_n)=\e[\mathds{1}_{\Lambda_n}M_n^\varphi]\;.
\end{equation*}
Then\\
\begin{enumerate}
\item $S_\infty$ is finite $Q^\varphi$-a.s. and satisfies for 
every $k\in\mathbb N$:
\begin{equation*}
Q^\varphi(S_\infty=k)=\varphi(k)\;.
\end{equation*}
\item Under $Q^\varphi$, the r.v. $T_\infty:=
\inf\left\lbrace n\geq0,\ X_n=S_\infty\right\rbrace $ 
(which is not a stopping time in general) is a.s. finite and\\
\begin{enumerate}
\item $(X_{n\wedge T_\infty},\ n\geq0)$ and 
$(S_\infty-X_{T_\infty+n},\ n\geq0)$ 
are two independent processes;\\
\item conditional on the r.v. $S_\infty$,
the process $(X_{n\wedge T_\infty},\ n\geq0)$ is a standard 
random walk stopped when it first hits the level $S_\infty$;\\
\item $(S_\infty-X_{T_\infty+n},\ n\geq0)$ is a 3-Bessel walk
started from $0$.\\
\end{enumerate}
\end{enumerate}
\item Put $R_n=2S_n-X_n$. Under $Q^\varphi$, $\left( R_n, n\geq0\right)$ 
is a 3-Bessel walk independent of $S_\infty$.
\end{enumerate}
\end{thm}

The 3-Bessel walk is the Markov chain $(R_n,n\geq0)$, with values
in  $\mathbb N=\{0,1,2,\ldots\}$, whose transition probabilities
from $x\geq0$ are given by
\begin{equation}
\pi(x,x+1)=\frac{x+2}{2x+2}\;; \qquad 
\pi(x,x-1)=\frac{x}{2x+2}\;.
\end{equation}

The 3-Bessel* walk is the Markov chain $( R^*_n,n\geq0)$,
valued in $\mathbb N^*=\{1,2,\ldots\}$, such that
$ R^*-1$ is a 3-Bessel walk.
So its transition probabilities from $x\geq1$ are
\begin{equation*}\label{transBes*}
\pi^*(x,x+1)=\frac{x+1}{2x}\;;\qquad
\pi^*(x,x-1)=\frac{x-1}{2x}\;.\end{equation*}\\
the 3-Bessel walk and the 3-Bessel* walk, will play a role 
in this work; they are identical up to a one-step space shift.\\
This result and those of \cite{D:1} can let think that the process of penalization gives rather intuitive results. Nevertheless, the following Theorems show that this intuition can be false and it is necessary to lead the calculations to their terms.

1) In the first section, $G$ is a function of the one-sided maximum till the last zero before $p$, i.e. :
$$G_p=\varphi(S_{g_p})$$
where $g_p:=\sup\left\lbrace k\leq p, X_k=0\right\rbrace$ and where $\varphi$ satisfies (\ref{conditions}) and  :
\begin{equation}\label{conditions2}
\sum_{k=0}^\infty k\varphi(k)<\infty .
\end{equation}
To study this penalisation we have to introduce $(\gamma_n,n\geq0)$ the number of $0$ before $n$ and we also recall that for all real $a$, $a^+:=\sup(a,0)$.
The result of this first section is summarized in the following statement :
\begin{thm}\label{thm5.1}
\begin{enumerate}
\item
\begin{enumerate} 
\item For all $n\geq0$ et all $\Lambda_n\in\mathcal F_n$:
\begin{equation}\label{f1.1}
\lim_{p\rightarrow\infty}\frac{E\left[ \mathds{1}_{\Lambda_n}\varphi\left( S_{g_p}\right) \right]}{ E\left[\varphi\left( S_{g_p}\right) \right]}=E\left[ \mathds{1}_{\Lambda_n}M_n\right], 
\end{equation}
where $M_n=\frac{1}{2}\varphi\left( S_{g_n}\right)\vert X_n\vert +\varphi\left( S_n\right) \left( S_n-X_n^+\right)+\phi\left( S_n\right)$.
\item Moreover, $\left(M_n,n\geq0 \right) $ is a positive martingale, not uniformly integrable.
\end{enumerate}
\item Let $Q$ be the probability on $\left(\Omega,\mathcal F_\infty \right) $, induces by:
\begin{equation*}
\forall n\geq0, \Lambda_n\in \mathcal F_n, Q\left(\Lambda_n \right):=E\left[ \mathds{1}_{\Lambda_n}M_n\right] .
\end{equation*}
Then under the probability $Q$:
\begin{enumerate}
\item let $g:=\sup\left\lbrace k\geq0, X_k=0\right\rbrace$. Then $Q\left(0\leq g<\infty \right)=1$.
\item $Q\left( S_\infty=\infty\right)=\frac{1}{2}$ and, conditionally on $S_\infty<\infty$, $\varphi$ is the density of $S_\infty$.
\item $\left( S_g,\gamma_g\right)$ admits as density:
$$f_{_{ \gamma_g,S_g }}(a,k):=\left\lbrace
  \begin{array}{cl}
    \left( \frac{1}{2}\right)^{a}\varphi(0)& \mbox{, for $k=0$}\\
\frac{1}{2}\left\lbrace\left(1-\frac{1}{2(k+1)} \right)^{^{a-1}}-\left(1-\frac{1}{2k} \right)^{^{a-1}} \right\rbrace \varphi(k)&\mbox{, otherwise.} 
  \end{array}
\right.$$
In particular, $S_g$ admits $\varphi$ as density.
\end{enumerate}
\item Under $Q$:
\begin{enumerate}
\item $\left( X_n,n\leq g\right)$ and $\left( X_n,n>g\right)$ are two independent processes.
\item With probability $\frac{1}{2}$, $\left( X_{g+n},n\geq0\right)$ (respectively to $\left( -X_{g+n},n\geq0\right)$) is a 3-Bessel* walk.
\item Conditionally on $\gamma_g=a$ and $S_g=b$, the process $\left( X_n,n\leq g\right)$ is a symmetric random walk stopped in $\tau_a$ and conditionally on $S_{\tau_a}=b$.
\end{enumerate}
\end{enumerate}
\end{thm}
2)  In the second section, $G_p=\varphi\left(S_{d_p}\right)$ where $d_p:=\inf\left\lbrace k\geq p,X_k=0\right\rbrace$ the first zero after $p$ and  $\varphi$  satisfies (\ref{conditions}) and (\ref{conditions2}). Let $f:\mathbb N\times\mathbb Z\rightarrow\mathbb R^+$ such that:
\begin{equation*}
f(b,a):=\mathds{1}_{b=0}\varphi(0)+\mathds{1}_{b\neq0}\left[ \varphi(b)\left(1-\frac{a^+}{b}\right)+a^+\sum_{k=b}^\infty\frac{\varphi(k)}{k(k+1)} \right] .
\end{equation*}
The main result of this section is :
\begin{thm}\label{thm6.1}
For all $n\geq0$ and all $\Lambda_n\in\mathcal F_n$:
\begin{equation*}
\lim_{p\rightarrow\infty}\frac {E\left[ \mathds{1}_{\Lambda_n}\varphi\left(S_{d_p} \right) \right] }{E\left[\varphi\left(S_{d_p} \right) \right] }=\lim_{p\rightarrow\infty}\frac {E\left[ \mathds{1}_{\Lambda_n}f\left( S_p,X_p\right) \right] }{E\left[ f\left( S_p,X_p\right)\right] }=\lim_{p\rightarrow\infty}\frac {E\left[ \mathds{1}_{\Lambda_n}\varphi\left(S_{p} \right) \right] }{E\left[\varphi\left(S_{p} \right) \right] }=E\left[\mathds{1}_{\Lambda_n}M_n^\varphi\right] ,
\end{equation*}
where $M_n^\varphi:=\varphi\left( S_n\right)\left( S_n-X_n\right) +1-\phi\left( S_n\right)$.\\
We remark that we obtain the same martingale as the one obtained for the penalisation by a function of the maximum( cf \cite{D:1}), i.e. that the penalisation by $\varphi\left( S_{d_p}\right)$ is the same as the penalisation by $\varphi\left( S_p\right) $.
\end{thm}
3) In the third section, $\varphi$ has to satisfy a stronger integrability condition :\begin{equation}
\sum_{k\geq0}k^2\varphi(k)<\infty
\end{equation}
and $G_p=\varphi\left(S^*_p\right)$ where $S^*_p=\sup_{n\leq p}\vert X_n\vert$, the bilateral supremum of $(X_n)_{n\geq0}$.
The result of this alinea is :
\begin{thm}\label{thm8.1}
\begin{enumerate}
\item 
\begin{enumerate}\item  Let $n\in\mathbb N$, $\Lambda_n\in\mathcal F_n$:
\begin{equation*}
\lim_{p\rightarrow\infty}\frac{E\left[ \mathds{1}_{\Lambda_n}\varphi\left( S^*_{g_p}\right) \right]}{ E\left[\varphi\left( S^*_{g_p}\right) \right]}=E\left[ \mathds{1}_{\Lambda_n}M^*_n\right] ,
\end{equation*}
where $M_n^*=\varphi\left( S^*_{g_n}\right)\vert X_n\vert+\varphi\left( S_n^*\right) \left(S_n^*-\vert X_n\vert \right)+\phi\left( S_n^*\right)$.
\item Moreover, $\left( M_n^*,n\geq0\right) $ is a non uniformly,  integrable positive $\mathcal F_n$ martingale .
\end{enumerate}
\item  Let $Q^*$ be the probability on $\left( \Omega,\mathcal F_\infty\right)$ induced by:
\begin{equation*}
\forall n\in\mathbb N,\Lambda_n\in\mathcal F_n:\quad Q^*\left( \Lambda_n\right):=E\left[\mathds{1}_{\Lambda_n} M^*_n\right] .
\end{equation*}
So, under $Q^*$:
\begin{enumerate}
\item Let $g:=\sup\left\lbrace k\geq0,X_n=0\right\rbrace $. then $g$ is finite and $S_\infty=\infty$ a.s.
\item The law of the couple $\left( S_g,\gamma_g\right)$ is:
$$f_{\gamma_g,S_g}\left( a,k\right)=\left\lbrace
\begin{array} {cl}
\varphi(0)&\mbox{, if $a=k=0$}\\
\left\lbrace\left(1-\frac{1}{k+1} \right)^a -\left(1-\frac{1}{k} \right) ^a \right\rbrace \varphi(k)&\mbox{, for $a\geq 0$, $ k>0$.}
\end{array}
\right.
$$
We deduced that $\varphi$ is the density of $S_g$.
\end{enumerate}
\item Under $Q^*$:
\begin{enumerate}
\item $\left( X_n,n\leq g\right)$ and $\left( X_n,n>g\right)$ are two independent processes.
\item With probability $\frac{1}{2}$, $\left( X_{g+n},n\geq0\right)$ (respectively $\left( -X_{g+n},n\geq0\right) $) is a three dimensional Bessel* walk.\\
\item Conditionally on $\gamma_g=a$ and $S^*_g=b$, $\left( X_n,n\leq g\right)$ is a symmetric random walk stopped in $\tau_a$ and conditionally on $S^*_{\tau_a}=b$.
\end{enumerate}
\end{enumerate}
\end{thm}

4) Finally, in order to be comprehensive, we fix an integer $a>0$ and consider the penalisation functional : 
$$G_p=\mathds{1}_{\left\lbrace S^*_p<a  \right\rbrace}.$$
We obtain :
\begin{thm}\label{thm7.1}
\begin{enumerate}
\item For each $n\geq 0$ and each $\Lambda_n\in\mathcal F_n$ :
\begin{equation*}
\lim_{p\rightarrow\infty}\frac{\e\left[\mathds{1}_{\left\lbrace\Lambda_n,\,S^*_p<a\right\rbrace}\right]}{\e\left[\mathds{1}_{\left\lbrace S^*_p<a\right\rbrace}\right]}:=\e\left[\mathds{1}_{\left\lbrace\Lambda_n,S^*_n<a\right\rbrace}M_n\right],
\end{equation*}
where $M_n:=\mathds{1}_{\left\lbrace\Lambda_n,S^*_n<a\right\rbrace}\left(\cos\left({\frac{\pi}{2a}}\right)\right)^{-n}\sin\left( \frac{\pi(a-X_n)}{2a}\right)$ is a positive martingale non uniformly integrable.
\item Let us define a new probability $Q$ on $(\Omega, \mathcal F_\infty)$ characterized by :
\begin{equation*}
\forall n\in\mathbb N,\,\forall \Lambda_n\in\mathcal F_n,\, Q\left( \Lambda_n\right):=\e\left[\Lambda_nM_n\right].
\end{equation*}
Under this new probability $Q$, $(X_n,n\geq0)$ has the following transition probabilities for $-b+1\leq k\leq a-1$:
\begin{eqnarray*}
Q\left(X_{n+1}=k+1\vert X_n=k\right)&=&\frac{\sin\left(\frac{a-k-1}{2a}\pi\right)}{2\cos\left(\frac{\pi}{2a}\right)\sin\left(\frac{a-k}{2a}\pi\right)},\\
Q\left(X_{n+1}=k-1\vert X_n=k\right)&=&\frac{\sin\left(\frac{a-k+1}{2a}\pi\right)}{2\cos\left(\frac{\pi}{2a}\right)\sin\left(\frac{a-k}{2a}\pi\right)}.
\end{eqnarray*}
\end{enumerate}
\end{thm}

\section{Penalisation by a function of $\mathbf{S_{g_p}}$, proof of Theorem \ref{thm5.1}}
1) To establish the first point of the Theorem (formula \ref{f1.1}), we need the following lemma :
\begin{lem}\label{newlem}
$\forall x\in\mathbb N,\, \forall a\in\mathbb Z\backslash ]x,+\infty[$ :
\begin{equation*}
\frac{\p_a\left( S_{g_n}=x,T_0<n\right)}{2\p_a(S_n=0)}
\end{equation*}
is bounded above by 1 for all $n\geq 0$ and tends to 1 when $n\rightarrow\infty$.
\end{lem}
\begin{proof}
We have to see that :
\begin{equation*}
\p_a\left( S_{g_n}<x,T_0<n\right)=\p_a\left(T_0<n,g_n<T_x\right)=\p_a\left(T_0<n<T_x+T_0\circ \theta_{T_x}\right),
\end{equation*}
where $\left\lbrace \theta_n\right\rbrace_n$ denote the family of shifts operators. We split in two cases according to the sign of $a$.\\
First, for  $a\leq 0$, according to the Desir\'e Andr\'e's principle :
$$\p_a\left( S_{g_n}<x,T_0<n\right)=\p_a\left(0\leq S_n<2x\right)=\p\left( \vert a\vert \leq S_n<2x+\vert a\vert  \right). $$
Which implies : 
\begin{eqnarray*}
\p_a\left(S_{g_n}=x,T_0<n\right)&=&\p\left( \vert a\vert \leq S_n<2x+2+\vert a\vert  \right) -\p\left( \vert a\vert \leq S_n<2x+\vert a\vert  \right) \\
&=&\p\left( S_n=2x+\vert a\vert  \right) + \p\left( S_n=2x+1+\vert a\vert  \right) .
\end{eqnarray*}
And for $a>0$ :
\begin{multline*}
\p_a\left( S_{g_n}<x,T_0<n\right)=\p_a\left(n<T_x+T_0\circ\theta_{T_x}, T_0<T_x \right)
-\p_a\left(n<T_x+T_0\circ\theta_{T_x},n\leq T_0<T_x\right)\\
=\p_a\left(n<T_x+T_0\circ\theta_{T_x}, T_0<T_x \right)-\p_a\left(n\leq T_0<T_x\right)\\
=\p_a\left(n<T_x+T_0\circ\theta_{T_x} \right)-\p_a\left(n<T_x+T_0\circ\theta_{T_x}, T_x<T_0 \right)
-\p_a\left(n\leq T_0<T_x\right)\\
=\p_a\left(n<T_x+T_0\circ\theta_{T_x} \right)-\p_a\left(n<T_0, T_x<T_0 \right)
-\p_a\left(n\leq T_0<T_x\right)\\
=\p_a\left(n<T_x+T_0\circ\theta_{T_x} \right)-\p_a\left(n<T_0\right)+\p_a\left(n<T_0<T_x\right)
-\p_a\left(n\leq T_0<T_x\right)\\
=\p_a\left(n<T_{2x} \right)-\p_a\left(n<T_0\right)-\p_a\left(n=T_0<T_x\right).
\end{multline*}
And consequently :
\begin{eqnarray*}
\p_a\left( S_{g_n}=x,T_0<n\right)&=&\p_a\left(2x\leq S_n<2x+2\right)-\p_a\left(n=T_0<T_{x+1}\right)+\p_a\left(n=T_0<T_x\right)\\
&=&\p_a\left(2x\leq S_n<2x+2\right)-\p_a\left(n=T_0,S_n=x\right)\\
&\leq&\p\left(2x-a\leq S_n<2x+2-a\right).
\end{eqnarray*}
In the ratio $\frac{\p(S_n=x)}{\p(S_n=0)}$, the denominator is bounded below by 
$\p(X_1=\ldots=X_n=-1)=2^{-n}$; so~it does not vanish.
Observe that, for even $n$ and even $k\geq2$,
$${\p(S_n=k{-}1)\over\p(S_n=0)}
={\p(S_n=k)\over\p(S_n=0)}={p_{n,k}\over p_{n,0}}
=\Bigl({n{-}k{+}2\over n{+}2}\Bigr)\ 
\Bigl({n{-}k{+}4\over n{+}4}\Bigr)\>\cdots\>
\Bigl({n\over n{+}k}\Bigr)\;;$$
and for odd $n$ and odd $k\geq 1$,
$${\p(S_n=k{-}1)\over\p(S_n=0)}
={\p(S_n=k)\over\p(S_n=0)}={p_{n,k}\over p_{n,1}}
=\Bigl({n{-}k{+}2\over n{+}1}\Bigr)\ 
\Bigl({n{-}k{+}4\over n{+}3}\Bigr)\>\cdots\>
\Bigl({n{+}1\over n{+}k}\Bigr)\;.$$
Clearly, these products are smaller than $1$
and tend to $1$ when $n$ goes to~infinity.
This proves the first point of the lemma. Obviously, when $a\leq 0$, the ratio tends to 1 when $n$ goes to infinity. 
In the other case, it appears clearly that $\p_a\left(n=T_0,S_n=x\right)\leq \p_a\left(I_n=0,S_n=x\right)$ tends to zero faster then the quantity $\p(S_n=0)$ (we have explicitely the expression of  $ \p_a\left(I_n=0,S_n=x\right)$ a little bit further in this paper).
\end{proof}
\begin{rem}\label{newrem}
Remark that we have also proved that for each $k\geq 0$ the ratio : $$\frac{\p(S_p=k)}{\p(S_p=0)}$$
 is bounded above by 1 and tends to 1 when $p\rightarrow+\infty$.
\end{rem}
\begin{lem}\label{newlem2}
For all $x\in\mathbb Z$ and $a\in\mathbb Z\backslash]x,+\infty[$ :
\begin{equation*}
\frac{\e_a\left[ \varphi\left(x\vee S_{g_n} \right)\mathds{1}_{T_0<n}\right]}{2\p(S_n=0)}
\end{equation*}
is bounded above by $(x-a^+)\varphi(x)+\phi(x) $ and tends to $(x-a^+)\varphi(x)+\phi(x)$ when $n\rightarrow\infty$.
\end{lem}
\begin{proof}
Write :
\begin{eqnarray*}
\frac{\e_a\left[\varphi\left(x\vee S_{g_n}\right)  \mathds{1}_{T_0<n}\right]}{2\p(S_n=0)}&=&\frac{\p_{a}\left(S_{g_n}<x,T_0<n\right)}{2\p(S_n=0)}\varphi(x)\\
&+&\sum_{k=x}^{\infty}\frac{\p_{a}\left(S_{g_n}=k,T_0<n\right)}{2\p(S_n=0)}\varphi(k).
\end{eqnarray*}
By lemma \ref{newlem}, this sum is bounded above by $(x-a^+)\varphi(x)+\phi(x) $ and tends to this value by dominated convergence.
\end{proof}
To prove point 1.a, we split:
\begin{eqnarray*}
\e\left[ \varphi\left( S_{g_p}\right)\mid\mathcal F_n\right]=\e\left[ \varphi\left( S_{g_p}\right)\mathds{1}_{g_p<n}\mid\mathcal F_n\right]+\e\left[ \varphi\left( S_{g_p}\right)\mathds{1}_{g_p\geq n}\mid\mathcal F_n\right]
:=(1)+(2).
\end{eqnarray*}
a) As $0\leq n\leq p$, we can write $(\tilde X_k,{k\geq 0}):=\left( X_{n+k}-X_n,{k\geq 0}\right)$, a standard random walk independent of $\mathcal F_n$. We denote by $\tilde T_a$ and $\tilde S$, the hitting time of the level $a$ and the supremum associated to $\tilde X$. Obviously on $\left\lbrace g_p<n \right\rbrace$, $\left\lbrace g_p=g_n\right\rbrace$. Hence :
$$(1)=\varphi(S_{g_n})\tilde \p_{X_n}(\tilde T_{0}>p-n)=\varphi(S_{g_n})\tilde\p\left(\tilde S_{p-n}\leq \vert X_n\vert \right), $$
where $\tilde \p$ only integrates over $\tilde S$, $X_n$ being kept fixed. 
Eventually, according to remark \ref{newrem} :
$$\frac{\e\left[\mathds{1}_{\Lambda_n,g_p<n}\varphi\left(S_{g_p}\right)\right]}{\p(S_{p-n}=0)}=\frac{\e\left[\mathds{1}_{\Lambda_n}\varphi\left(S_{g_n}\right)\tilde\p\left(\tilde S_{p-n}\leq \vert X_n\vert \right)  \right]}{\p(S_{p-n}=0)} $$
is bounded above by $\e\left[\mathds{1}_{\Lambda_n}\varphi(S_{g_n})\vert X_n\vert\right]$ which is integrable and tends to $\left[\mathds{1}_{\Lambda_n}\varphi(S_{g_n})\vert X_n\vert \right]$ when $p$ goes to $\infty$. 

b) We now study the behaviour of $(2)$. We use the same notations as before, adding for all $p\geq 0$, $\tilde g_p$ the last zero before $p$ associated to $\tilde X$. Hence :
\begin{equation*}
(2)=\frac{\e\left[\e\left[ \varphi\left( S_{n}\vee S_{\left[ n,g_p\right] }\right)\mathds{1}_{g_p\geq n}\mid\mathcal F_n\right]\right]}{2\p(S_{p-n}=0)}
=\e\left[\frac{\tilde \e_{X_n}\left[ \varphi\left( S_{n}\vee \tilde S_{\tilde g_{p-n}} \right)\mathds{1}_{\tilde T_{0}\leq p-n}\right]}{2\p(S_{p-n}=0)}\right],
\end{equation*}
where $\tilde \e$ integrates on $\tilde S$, $\tilde g$ and $\tilde T_0$,  the variables $S_n$ and $X_n$ being kept fixed.  
When $p$ tends to infinity, Lemma \ref{newlem} says that the ratio in the right hand side tends to $(S_n-X_n^+)\varphi(S_n)+\phi(S_n)$ and is dominated by the same quantity, which is integrable. As a result :
$$\frac{\e\left[\mathds{1}_{\Lambda_n}\varphi(S_{g_p})\right]}{\p(S_{p-n}=0)}
\left\lbrace 
\begin{aligned}
 \mbox{is bounded above by $\e\left[\mathds{1}_{\Lambda_n}M_n\right]$ forall $p\geq n$.}\\
 \mbox{and tends to  $\e\left[\mathds{1}_{\Lambda_n}M_n\right]$ when $p\rightarrow\infty$.}
\end{aligned}
\right.
$$
Taking in particular $\Lambda_n=\Omega$, one also has 
$$\frac{\e\left[\varphi\left(S_{g_p}\right)\right]}{2\p(S_{p-n}=0)} \underset{p\rightarrow\infty}{\rightarrow}\e\left[M_n\right]=1,$$
and to establish point 1 of Theorem \ref{thm5.1}, it suffices to take the ratio of these two limits.\\
ii) Let us prove now that $\left(M_n,n\geq0 \right) $ is a $(\mathcal F_n)$-martingale under $\p$.
For typographical simplicity, we write :
$$M_n=\mathcal A_n+\mathcal B_n $$
where  $\mathcal A_n:=\varphi(S_{n})(S_{n}-X_{n}^+)+\sum_{k=S_n}^\infty\varphi(k)$ and $\mathcal B_n:=\frac{1}{2}\varphi(S_{g_{n}})\vert X_{n}\vert$.\\
We suppose that $n>0$, the case $n=0$ being trivial.\\
On $\left\lbrace X_n\geq 1\right\rbrace$, $\mathcal A_{n+1}$ is in fact the martingale found in the Theorem \ref{thmmax}, then conditional on $\mathcal F_n$, this quantity is equal to $\mathcal A_n$.\\
On $\left\lbrace X_n\leq -1\right\rbrace$, $S_{n+1}=S_n=S_{g_n}$ and $X_{n+1}^+=X_n^+=0$, obviously on this event  $\mathcal A_{n+1}=\mathcal A_n$.\\
Eventually, on $\left\lbrace X_n=0\right\rbrace$, as $\vert X_{n+1}\vert=1$, we have $S_{n+1}=S_n$. So, summing on the possible values of $X_{n+1}$, 1 and $-1$, it is easy to check that :
\begin{equation}\label{martA}
\e[\mathds{1}_{X_{n}=0}\mathcal A_{n+1}\vert \mathcal F_n]=\mathds{1}_{X_{n}=0}\left(\mathcal A_{n}-\frac{1}{2}\varphi(S_{g_n})\right).
\end{equation}
It just remains the quantity $\mathcal B_n:=\frac{1}{2}\varphi(S_{g_{n}})\vert X_{n}\vert$.\\
On $\left\lbrace\vert X_n\vert\geq 2\right\rbrace$, $S_{g_{n+1}}=S_{g_n}$ and as the function $x\rightarrow\vert x\vert$ is harmonic for the symmetric random walk, except in $0$, consequently $\e[\mathds{1}_{\vert X_n\vert \geq 2}\mathcal B_{n+1}\vert\mathcal F_n]=\mathds{1}_{\vert X_n\vert\geq 2}\mathcal B_n$.\\
On $\left\lbrace\vert X_n\vert =1\right\rbrace$, either $\vert X_{n+1}\vert=2$ and in this case $S_{g_{n+1}}=S_{g_n}$ implies $\mathcal B_{n+1}=\varphi( S_{g_n})$, either $\vert X_{n+1}\vert=0$ and in this case $\mathcal B_{n+1}=0$. Then, immediately we have $\e[\mathds{1}_{\vert X_n\vert =1}\mathcal B_{n+1}\vert \mathcal F_n]=\mathds{1}_{\vert X_n\vert=1}\frac{1}{2}\varphi(S_{g_n})=\mathds{1}_{\vert X_n\vert=1}\mathcal B_n$.\\
At last, on $\left\lbrace X_n=0\right\rbrace$, $S_{g_{n+1}}=S_{g_n}$ and consequently $\mathcal B_{n+1}=\frac{1}{2}\varphi(S_{g_n})$. So according to (\ref{martA}) :
\begin{equation*}
\e[\mathds{1}_{X_n=0}M_{n+1}\vert \mathcal F_n]=\e[\mathds{1}_{X_n=0}(\mathcal A_{n+1}+\mathcal B_{n+1})\vert \mathcal F_n]=\mathds{1}_{X_n=0}\mathcal A_{n}=\mathds{1}_{X_n=0}M_{n}.
\end{equation*}\\
2) For $p$ and $n$ in $\mathbb N$, the event $\left\lbrace S_n>p\right\rbrace$ is equal to $\left\lbrace T_p<n\right\rbrace$. Using the definition of $Q$ and the Doob's stopping Theorem :
\begin{equation*}
Q\left( S_n>p\right)=Q\left( T_p<n\right)=\e[\mathds{1}_{T_p<n}M_{T_p}]=\e\left[\mathds{1}_{T_p<n}\left\lbrace \frac{1}{2}\varphi\left( S_{g_{T_p}}\right)p+\phi\left(p \right)   \right\rbrace \right].
\end{equation*}
Moreover according to \cite{LG:1} p.457-458,
under $\p$, $S_{g_{T_p}}$ is a uniformly distributed random variable on $\left\lbrace0,...,p-1 \right\rbrace$.
As  $n\rightarrow\infty$, the Lebesgue Theorem permits us to write:
\begin{eqnarray*}
Q\left( S_\infty>p\right)&=&\lim_{n\rightarrow\infty}Q\left( S_n>p\right)
=\e\left[\left\lbrace \frac{1}{2}\varphi\left( S_{g_{T_p}}\right)p+\phi\left(p \right)   \right\rbrace \right]\\
&=&\frac{p}{2}\e\left[\varphi\left( S_{g_{T_p}}\right)\right]+\phi\left(p \right)
=\frac{p}{2}\sum_{k=0}^{p-1}\frac{1}{p}\varphi(k)+\phi(p)
=\frac{1}{2}\sum_{k=0}^{p-1}\varphi(k)+\phi(p).
\end{eqnarray*}
Consequently 
$
Q\left( S_\infty=\infty\right)=\lim_{p\rightarrow\infty}Q\left( S_\infty>p\right)
=\frac{1}{2}
$
and the half of the point 2.b is proved.\\
In order to prove point 2.a,
we need, for $a>0$, the law under $\p$ of $S_{d_a}$ conditionally on $\mathcal F_p$.
\begin{lem}\label{lem5.5}
Let $k\geq a>0$, then:
\begin{equation*}
P_a\left(S_{T_0}=k \right)=\frac{a}{k(k+1)}.
\end{equation*}
\end{lem}
\begin{proof}
A direct use of the stopping Theorem to the martingale $(X_n,n\geq 0)$ and the stopping time $T_0\wedge T_k$ gives us : 
\begin{equation*}
\p_a\left( T_0>T_k\right)=\frac{a}{k}.
\end{equation*}
We just have to remark that $\p_a\left(S_{T_0}=k \right)=\p_a\left(T_k<T_0\right)-\p_a\left(T_{k+1}<T_0\right)$ to achieve the proof.
\end{proof}
\begin{lem}\label{lem5.6}
Let $\psi:\mathbb N\rightarrow\mathbb R^+$ be an integrable function. Then :
\begin{equation*}
\e\left[\psi\left( S_{d_p}\right)\mid\mathcal F_p  \right]=\mathds{1}_{S_p=0}\psi(0)+ \mathds{1}_{S_p\neq0}\left\lbrace  \psi\left( S_p\right)\left( 1-\frac{X_p^+}{S_p}\right)+X_p^+\sum_{k\geq S_p}\frac{\psi(k)}{k(k+1)} \right\rbrace.
\end{equation*}
\end{lem}
\begin{proof}
We easily obtain :
\begin{equation*}
\e\left[\psi\left( S_{d_p}\right)  \mid\mathcal F_p\right]
=\mathds{1}_{S_p=0}\psi\left( 0\right)+\e\left[\mathds{1}_{S_p\neq0}\psi\left( S_p\vee S_{\left[p,d_p\right]}\right)\mid\mathcal F_p  \right].
\end{equation*}
If $X_p\leq 0$, then $X_k\leq 0$ for all $p\leq k\leq d_p$ and consequently $X_p^+=S_{[p,d_p]}=0$.
As a result, on $\left\lbrace X_p\leq 0, S_p\neq 0\right\rbrace$ :
\begin{equation*}
\psi\left( S_p\vee S_{\left[p,d_p\right]}\right)=\psi(S_p)=\psi(S_p)\left(1-\frac{X_p^+}{S_p}\right)+X_p^+\sum_{k\geq S_p}\frac{\psi(k)}{k(k+1)}.
\end{equation*}
Let $( \tilde X_q:=X_{q+p}\geq0)$, a random walk starting from $X_p$ and independent of $\mathcal F_p$ and we denote by $\tilde S$ and $\tilde T_0$ respectively the supremum and the hitting time of 0 associated to $\tilde X$.
Then $S_{[p,d_p]}=\tilde S_{\tilde T_0}$. 
In the following calculus, $\tilde \e$ only integrates $\tilde S_{\tilde T_0}$, $X_p$ and $S_p$ being kept fixed. Consequently, on $\left\lbrace X_p>0, S_p\neq 0\right\rbrace$, according to lemma \ref{lem5.5}, $\e\left[\psi\left( S_{d_p}\right)  \mid\mathcal F_p\right]$ equals to:
\begin{eqnarray*}
\tilde \e_{X_p}\left[\mathds\psi\left( S_p\vee\tilde S_{\tilde T_0}  \right)\right]
&=&\sum_{k=X_p}^{S_p-1}\tilde \p_{X_p}\left(\tilde S_{\tilde T_0}=k \right) \psi\left( S_p\right)+\sum_{k\geq S_p}^{\infty}\tilde \p_{X_p}\left(\tilde S_{\tilde T_0}=k \right) \psi\left(  k  \right)\\
&=& \sum_{k=X_p}^{S_p-1}\frac{X_p}{k(k+1)}\psi\left( S_p\right)+\sum_{k\geq S_p}^{\infty} \frac{X_p}{k(k+1)}\psi\left(  k  \right) \\
&=& \left( 1-\frac{X_p}{S_p}\right)\psi\left( S_p\right)+\sum_{k\geq S_p}^{\infty} \frac{X_p}{k(k+1)}\psi\left(  k  \right).
\end{eqnarray*}
\end{proof}
Fixing $a>0$, according to Doob's stopping  Theorem:
\begin{eqnarray*}
Q\left( g_p>a\right)=Q\left( d_a<p\right)
=\e\left[\mathds{1}_{d_a<p}M_{d_a} \right]
= \e\left[\mathds{1}_{d_a<p}\left\lbrace  \varphi\left( S_{d_a}\right) S_{d_a}+\phi\left( S_{d_a}\right) \right\rbrace \right].
\end{eqnarray*}
The events $\left\lbrace  g_p>a\right\rbrace$ form an increasing sequence with limit $\left\lbrace g>a\right\rbrace$. Hence :
\begin{eqnarray*}
Q\left( g>a\right)=\lim_{p\rightarrow\infty}Q\left( g_p>a\right)
=\e\left[ \varphi\left( S_{d_a}\right)  S_{d_a}+\phi\left( S_{d_a}\right) \right]
=\e\left[ \varphi\left( S_{d_a}\right)  S_{d_a}\right] +\e\left[ \phi\left( S_{d_a}\right) \right].
\end{eqnarray*}
To achieve the proof of the point 2.a,  we have to prove that each term tends to zero as $a\rightarrow\infty$.
According to the Lebesgue Theorem 
$\e\left[ \phi\left( S_{d_a}\right) \right]\underset{a\rightarrow\infty}{\rightarrow}0$.
We use lemma \ref{lem5.6} with $\psi(x):=x\varphi(x)$:
\begin{eqnarray*}
\e\left[ S_{d_a}\varphi\left( S_{d_a}\right) \right]
&=& \e\left[\mathds{1}_{S_a\neq0}\varphi\left( S_{a}\right)\left( S_a-X_a^+\right) +X_a^+\sum_{k\geq S_a}\frac{\varphi(k)}{k+1}\right]\\
&\leq& \e\left[\mathds{1}_{S_a\neq0}\varphi\left( S_{a}\right)S_a  +\frac{X_a^+}{S_a+1}\sum_{k\geq S_a}\varphi(k)\right]\\
&\leq& \e\left[\varphi\left( S_{a}\right)S_a +\phi\left( S_a\right) \right]
= \e\left[\varphi\left( S_{a}\right)S_a\right]  +\e\left[ \phi\left( S_a\right) \right].
\end{eqnarray*}
On the one hand, according to the Lebesgue Theorem 
$\e\left[ \phi\left( S_a\right) \right]\underset{a\rightarrow\infty}{\rightarrow}0$.\\
On the other hand $\e\left[\varphi\left( S_{a}\right)S_a\right]=\sum_{k=0}^\infty \p\left(S_a=k \right) \varphi(k)k$
is bounded above by  $\sum_{k=0}^\infty\varphi(k)k<\infty$ and $\p(S_a=k)$ tends to 0 when $a\rightarrow\infty$. Again, according to the Lebesgue Theorem,  $\e\left[\varphi\left( S_{a}\right)S_a\right]$ tends to 0 when $a\rightarrow\infty$.\\
3) First of all, let us establish preliminary results and remind that
$\gamma_{n}:=\sum_{k=0}^n\mathds{1}_{X_k=0}$ is the number of visits to zero up to time $n$ and  denote $ \tau_a:=\inf\left\lbrace p\geq0,\gamma_n=a\right\rbrace$.
\begin{lem}\label{lem5.7}
For all $c>0$ and $a\geq 1$:
\begin{equation*}
\p\left( S_{\tau_a}=c\right)=\left\lbrace\begin{array}{cc}
\left( \frac{1}{2}\right)^{a-1}&\mbox{, if $c=0$}\\
\left( 1-\frac{1}{2(c+1)}\right)^{a-1}-\left( 1-\frac{1}{2c}\right)^{a-1}&\mbox{, otherwise.}
\end{array}\right.
\end{equation*}
\end{lem}
\begin{lem}\label{lem5.8}
For all $n\geq0$:
\begin{equation*}
\sum_{k=n}^\infty \frac{1}{k(k+1)}\left[k\varphi(k)+\phi(k) \right]=\frac{1}{n}\left( 1-\phi(n)\right) .
\end{equation*}
\end{lem}
\begin{proof}[\textit{Proof of lemma \ref{lem5.7}}]
This is obvious for $a=1$ so let us suppose that $a\geq 2$.
If $c=0$, we have, with an obvious recurrence :
\begin{eqnarray*}
\p\left( S_{\tau_a}=0\right)=\p\left(X_1=-1 \right) \p\left( S_{\tau_a}=0\mid X_1=-1\right)
=\frac{1}{2}\p\left( S_{\tau_{a-1}}=0\right)
=\left(\frac{1}{2}\right)^{a-1}.
\end{eqnarray*}
Now suppose that $c>0$.
With those notations, using the strong Markov property and an obvious recurrence:
\begin{multline*}
\p\left( S_{\tau_a}<c\right)=\p \left( \tau_a<T_c\right)=\p(\tau_a<T_c\vert \tau_2<T_c)\p(\tau_2<T_c)=\p(\tau_{a-1}<T_c)\p(\tau_2<T_c)\\
=\p(\tau_2<T_c)^{a-1}=\left[\frac{1}{2}(\p_1(\tau_1<T_c)+\p_{-1}(\tau_1<T_c))\right]^{a-1}
=\left[\frac{1}{2}(\p_1(\tau_1<T_c)+1)\right]^{a-1}.
\end{multline*}
We have already seen  that $\p_1\left( T_c<T_0\right)=\frac{1}{c}$, then
$\p\left( S_{\tau_a}<c\right)=\left(1-\frac{1}{2c} \right)^{a-1}$.\\
We can note that the law of $\gamma_{T_c}$ is a geometric law of parameter $\frac{1}{2c}$.
Finally:
\begin{equation*}
\p\left( S_{\tau_a}=c\right)=\p\left( S_{\tau_a}<c\right)-\p\left( S_{\tau_a}<c+1\right)=\left(1-\frac{1}{2(c+1)} \right)^{a}-\left(1-\frac{1}{2c} \right)^{a}.
\end{equation*}
\end{proof}
\begin{proof}[\textit{Proof of lemma \ref{lem5.8}}]
We have :
\begin{eqnarray*}
\sum_{k=n}^\infty \frac{\phi(k)}{k(k+1)}&=&\sum_{k=n}^\infty\sum_{l}^\infty\frac{\varphi(l)}{k(k+1)}=\sum_{l=n}^\infty\sum_{k=n}^l \frac{\varphi(l)}{k(k+1)}=\sum_{l=n}^\infty\varphi(l)\left(\frac{1}{n}-\frac{1}{k+1}\right),
\end{eqnarray*}
hence :
\begin{equation*}
\sum_{k=n}^\infty \frac{1}{k(k+1)}\left[k\varphi(k)+\phi(k) \right]= \frac{1}{n} \sum_{k=n}^\infty\varphi(k).
\end{equation*}
\end{proof}
Let $F$ be a functional, $f_1$ and $f_2$ be two functions from $\mathbb N$ to $\mathbb R^+$.
\begin{eqnarray*}
\mathcal A &:=&\e^Q\left[ F\left( X_u,u\leq g\right)f_1\left(\gamma_{g} \right)f_2\left(S_g\right)\right]\\
&=&\sum_{a\geq1}   \e^Q\left[ F\left( X_u,u\leq \tau_a\right)f_1\left(\gamma_{\tau_a} \right)f_2\left(S_{\tau_a}\right)\mathds{1}_{\tau_{a}<\infty,\,\tau_{a+1}=\infty}\right]\\
&=&\sum_{a\geq1} \e^Q\left[ F\left( X_u,u\leq \tau_a\right)f_1\left(\gamma_{\tau_a} \right)f_2\left(S_{\tau_a}\right)(\mathds{1}_{\tau_{a}<\infty}-\mathds{1}_{\tau_{a+1}<\infty})\right]\\
&=&\sum_{a\geq1} \e\left[ F\left( X_u,u\leq \tau_a\right)f_1\left(\gamma_{\tau_a} \right)f_2\left(S_{\tau_a}\right)\mathds{1}_{\tau_{a}<\infty}M_{\tau_a}\right]\\
&-&\e\left[ F\left( X_u,u\leq \tau_a\right)f_1\left(\gamma_{\tau_a} \right)f_2\left(S_{\tau_a}\right)\mathds{1}_{\tau_{a+1}<\infty}M_{\tau_{a+1}}\right]\\
&=&\sum_{a\geq1} \e\left[ F\left( X_u,u\leq \tau_a\right)f_1\left(\gamma_{\tau_a} \right)f_2\left(S_{\tau_a}\right)\left( M_{\tau_a}-M_{\tau_{a+1}}\right) \right].
\end{eqnarray*}
Since:
\begin{equation*}
M_{\tau_a}-M_{\tau_{a+1}}
=\varphi\left( S_{\tau_a}\right)  S_{\tau_a}+\phi\left( S_{\tau_a}\right)-\varphi\left( S_{\tau_{a+1}}\right) S_{\tau_{a+1}}-\phi\left( S_{\tau_{a+1}}\right).
\end{equation*}
One can write $ S_{\tau_a+1}=S_{\tau_a}\vee\tilde S_{\tilde \tau_2}$ where $\tilde S$ and $\tilde \tau_2$ are respectively the unilateral supremum and the time of the return in 0 of the standard random walk $(X_{n+\tau_a},n\geq 0)$ which is independent of $\mathcal F_{\tau_a}$. Hence :
\begin{equation*}
M_{\tau_a}-M_{\tau_{a+1}}=\mathds{1}_{\tilde S_{\tilde \tau_2}>S_{\tau_a}}\left(\varphi(S_{\tau_a})S_{\tau_a}-\varphi(\tilde S_{\tilde \tau_2})\tilde S_{\tilde \tau_2}+\sum_{k=S_{\tau_a}}^{\tilde S_{\tilde \tau_2}-1}\varphi(k)\right) .
\end{equation*}
Then, we condition this quantity by $\mathcal F_{\tau_a}$ :
\begin{eqnarray*}
\e\left[M_{\tau_a}-M_{\tau_{a+1}}\vert \mathcal F_{\tau_a}\right]&=&\sum_{l>S_{\tau_a}}\p(\tilde S_{\tilde \tau_1}=l)(\varphi(S_{\tau_a})S_{\tau_a}-\varphi(l)l + \sum_{k=l}^{l-1}\varphi(k))\\
&=&\sum_{l>S_{\tau_a}}\frac{1}{2l(l+1)}(\varphi(S_{\tau_a})S_{\tau_a}-\varphi(l)l + \sum_{k=l}^{l-1}\varphi(k))\\
&=&\frac{\varphi(S_{\tau_a})S_{\tau_a}}{2(S_{\tau_a}+1)}-\sum_{l>S_{\tau_a}}\frac{\varphi(l)}{2(l+1)}+\sum_{k\geq S_{\tau_a}}\varphi(k)\sum_{l\geq k+1}\frac{1}{2l(l+1)}\\
&=&\frac{1}{2}\left(\frac{\varphi(S_{\tau_a})S_{\tau_a}}{S_{\tau_a}+1}-\sum_{l>S_{\tau_a}}\frac{\varphi(l)}{l+1}+\sum_{k\geq S_{\tau_a}}\frac{\varphi(k)}{k+1}\right)=\frac{\varphi(S_{\tau_a)}}{2}.
\end{eqnarray*}
Consequently:
$
\mathcal A
=\sum_{a\geq1} \frac{1}{2}\e\left[ F\left( X_u,u\leq \tau_a\right)f_1\left(a\right)f_2\left(S_{\tau_a}\right)\varphi\left( S_{\tau_a}\right)\right]
$ and
with $F\equiv1$ :
\begin{multline*}
\mathcal A
=\frac{1}{2}\sum_{a\geq1}\sum_{k\geq0}\p\left(S_{\tau_a}=k \right)f_1(a)f_2(k)\varphi(k)
=\frac{1}{2}\sum_{a\geq1}\left( \frac{1}{2}\right)^{a-1} f_1(a)f_2(0)\varphi(0)\\
+\frac{1}{2}\sum_{a\geq1}\sum_{k\geq1}\left\lbrace\left(1-\frac{1}{2(k+1)} \right)^{a-1}-\left(1-\frac{1}{2k} \right)^{a-1} \right\rbrace f_1(a)f_2(k)\varphi(k)
\end{multline*}
which gives us the density of $\left( \gamma_g,S_g \right)$.\\
Now, summing over $a$ we easily find that $\varphi$ is the density of $S_g$ under $Q$.\\
For proving 3.iii, we write the formula $\mathcal A$
in two different ways:
\begin{eqnarray*}
\mathcal A&=&\sum_{a\geq1}\sum_{k\geq0}f_{ \gamma_g,S_g }(a,k)\e_Q\left[ F\left( X_u,u\leq g\right)\mid S_g=k,\gamma_g=a\right]f_1\left(a \right)f_2\left(k\right)\\
&=&\frac{1}{2}\sum_{a\geq1}\sum_{k\geq0}f_1\left(a \right)f_2\left(k\right) \varphi\left( k\right)\p\left( S_{\tau_a}=k\right) \e\left[ F\left( X_u,u\leq \tau_a\right)\mid S_{\tau_a}=k \right]\\
\end{eqnarray*}
The formulas that we obtained for $Q\left(S_g=k,\gamma_g=a \right)$ and $\p\left( S_{\tau_a}=k\right)$ imply obviously that for all $k,a\geq0$:
\begin{equation*}
\e_Q\left[ F\left( X_u,u\leq g\right)\mid S_g=k,\gamma_g=a\right]=\e\left[ F\left( X_u,u\leq \tau_a\right)\mid S_{\tau_a}=k \right].
\end{equation*}
3.ii) The study of the process  $\left(X_n,\ n\geq0 \right) $ under $Q^{h^+,h^-}$
starts with the next three lemmas.
\begin{lem}\label{lem3}
Under $\p_1$  and conditional on the event 
$\left\lbrace T_p< T_0\right\rbrace $, the
process  $(X_n,\ 0\leq n\leq T_p)$ is a 3-Bessel* walk  started from 1 and stopped when it first hits the level $p$ (cf. \cite{LG:1}).
\end{lem}

For typographical simplicity, call 
 $T_{p,n}:=\inf\lbrace k> n,\ X_k=p\rbrace$ the
 time of the first visit to $p$ after $n$, and 
$\mathcal H_l:=\bigl\lbrace T_{p,\tau_l}<\tau_{l+1,\,X_{\tau_l+1}=1}\bigr\rbrace$, 
the event that the $l$-th excursion is positive 
and reaches level $p$.
\begin{lem}\label{lem4}
Under the law $Q$ and  conditional on 
the event $\mathcal H_l$, the
process $(X_{n+\tau_l},\ {1\leq n\leq T_{p,\tau_l}-\tau_l})$
 is a 3-Bessel* walk started from 1 and stopped when it first hits the level $p$.
\end{lem}
\begin{proof}
Let $G$ be a functional on $\mathbb Z^n$ :
\begin{eqnarray*}
\mathcal D&:=&Q\left[ G\left( X_{\tau_l+1},...,X_{\tau_l+n}\right)\mathds{1}_{n+\tau_l<T_{p,\tau_l}}\mid \mathcal H_l  \right]\\
&=&\frac{Q\left[ G\left( X_{\tau_l+1},...,X_{\tau_l+n} \right)\mathds{1}_{n+\tau_l<T_{p,\tau_l}<\tau_{l+1},X_{\tau_l+1}=1}\right]}{Q\left( \mathcal H_l\right)}\\
&=&\frac{\e\left[ G\left( X_{\tau_l+1},...,X_{\tau_l+n} \right)\mathds{1}_{n+\tau_l<T_{p,\tau_l}<\tau_{l+1},X_{\tau_l+1}=1}M_{\tau_{l+1}}\right]}{\e\left[ \mathds{1}_{\mathcal H_l}M_{\tau_{l+1}}\right]}.
\end{eqnarray*}
Obviously $M_{\tau_{l+1}}=\varphi(S_{\tau_{l+1}})S_{\tau_{l+1}}+\phi(S_{\tau_{l+1}})$ and conditionning by $\mathcal F_{T_{p,\tau_l}}$, we obtain $\e\left[ M_{\tau_{l+1}}\mid\mathcal F_{T_{p,\tau_l}}\right]=\e_p\left[\varphi(S_{T_0})S_{T_0}+\phi(S_{T_0})\right]$, a constant. Conditioning by $\mathcal F_{T_{p,\tau_{l}}}$ the denominator and numerator of $\mathcal D$ :
\begin{eqnarray*}
\mathcal D&=&\frac{\e\left[ G\left( X_{\tau_l+1},...,X_{\tau_l+n} \right)\mathds{1}_{n+\tau_l<T_{p,\tau_l}<\tau_{l+1},X_{\tau_l+1}=1}\e\left[M_{\tau_{l+1}}\mid \mathcal F_{T_{p,\tau_l}}\right]\right]}{\e\left[ \mathds{1}_{\mathcal H_l}\e\left[M_{\tau_{l+1}}\mid\mathcal F_{T_{p,\tau_l}}\right]\right]}\\
&=&\frac{\e\left[ G\left( X_{\tau_l+1},...,X_{\tau_l+n} \right)\mathds{1}_{n+\tau_l<T_{p,\tau_l}<\tau_{l+1},X_{\tau_l+1}=1}\right]}{\e\left[ \mathds{1}_{\mathcal H_l}\right]}
\end{eqnarray*}
Using the conditionning by $\mathcal F_{\tau_{l}+1}$ and the Markov property :
\begin{equation*}
\mathcal D=\frac{\e_1\left[ G\left( X_{0},...,X_{n-1} \right)\mathds{1}_{n-1<T_p<T_0}\right]}{\p_1\left( T_p<T_0\right)}=\e_1\left[G\left( X_{0},...,X_{n-1} \right)\mathds{1}_{n-1<T_p}\mid T_p<T_0\right].
\end{equation*}
\end{proof}
On the other part, according to \cite{LG:1} conditionally on $\left\lbrace T_p<T_0\right\rbrace $, the law of $\left( X_n , n< T_p\right) $  under $\p_1$ is the law of the 3-dimensional Bessel* walk. We deduce that, conditionally on  $\left\lbrace T_p<T_0\right\rbrace $ under $Q_1$, $\left( X_n , n< T_p\right) $ is a 3-dimensional Bessel* walk. Making  $p$ go to infinity, we obtain that under $Q_1$ conditionally on $\left\lbrace T_0=\infty\right\rbrace $, $\left( X_n , n\geq0\right) $ is a 3-dimensional Bessel* walk.\\
Obviously, by symmetry, under $Q_{-1}$, conditionally on $\left\lbrace T_0=\infty\right\rbrace $, $\left( -X_n , n\geq0\right) $ is a three dimensional Bessel* walk. We deduce that $\left(X_n,n\geq g \right)$ , is either a three dimensional Bessel* walk, either a reversed three dimensional Bessel* walk. It remains to know with what probability we have one or the other.\\
We have seen in 2.ii that $S_\infty$ under $Q$ was finished with probability $\frac{1}{2}$. As a three dimensional Bessel* walk goes to infinity in infinity, we deduce that  $\left(X_n,n\geq g \right)$ is one or the other walk with probability $\frac{1}{2}$.
\setcounter{equation}{0}

\section{Penalisation by a function of $\mathbf{S_{d_p}}$}
We've already seen according to lemma \ref{lem5.6} that $\e\left[\varphi\left( S_{d_p}\right)\mid\mathcal F_p  \right]=f\left( S_p,X_p\right) $. Moreover :
\begin{eqnarray*}
\e\left[ f\left( S_p,X_p\right)\mathds{1	}_{\Lambda_n}\right]&=&\e\left[ \varphi\left( S_p\right)\mathds{1	}_{\Lambda_n} \right] -\e\left[\mathds{1}_{S_p\neq0}\varphi\left( S_p\right)\frac{X_p^+}{S_p} \mathds{1	}_{\Lambda_n}\right]\\
&+&\e\left[\mathds{1}_{S_p\neq0} X_p^+\sum_{k\geq S_p}\frac{\varphi(k)}{k(k+1)} \mathds{1	}_{\Lambda_n}\right]\\
&=&(1)-(2)+(3).
\end{eqnarray*}
We already know (cf. \cite{D:1}) that  $\forall n\geq 0$ and $\Lambda_n\in\mathcal F_n$ :
\begin{equation}\label{limphi}
\frac{\e[\mathds1_{\Lambda_n}\,\varphi(S_p)]
}{\p(S_{p-n}=0)}\quad\left\{
\begin{array}{l}
\hbox{is bounded above by $\e[\mathds1_{\Lambda_n}\>M^\varphi_n]$
for all $p\ge n$}\\ \noalign{\vskip4pt}
\hbox{and tends to\quad}\e[\mathds1_{\Lambda_n}\>M^\varphi_n]
\hbox{\quad when }p\to\infty\;.
\end{array}\right.
\end{equation}
Then we just have to prove that :
$$\frac{\e\left[(f(S_p,X_p)-\varphi(S_p))\mathds{1}_{\Lambda_n}\right]}{\p(S_{p-n}=0)}$$
goes to 0 when $p\rightarrow\infty$.\\
In particular, if we take $\Lambda_n=\Omega$, we have :
$$\frac{\e\left[\varphi(S_{d_p})\mathds{1}_{\Lambda_n}\right]}{\p(S_{p-n}=0)}\underset{p\rightarrow\infty}{\rightarrow}\e\left[M_n^\varphi\right]=1 $$
and to establish Theorem \ref{thm6.1}, we take the ratio of the two limits.\\

ii) To study the behaviour of the last two terms, we need the following lemma:
\begin{lem}{\label{lem6.1}}
For $b\geq0$ and $a\leq b$ :
\begin{equation*}
\frac{\p\left( S_p=b,X_p=a\right)}{\p(S_p=0)}
\end{equation*}
is bounded above by 1 and tends to 0 when $p\rightarrow\infty$.
\end{lem}
\begin{proof}
Remark that $a$ and $p$ must have the same parity, otherwise\newline $\p\left( S_p=b,X_p=a\right)$ is equal to zero and the lemma is obvious.
According to remark \ref{newrem} :
$$\frac{\p(S_p=b,X_p=a)}{\p(S_p=0)}\leq  \frac{\p(S_p=b)}{\p(S_p=0)}\leq 1.$$
With these hypothesis, according to the Desir\'e  Andr\'e's reflexion principle:
\begin{multline*}
\p\left( S_p=b,X_p=a\right)=\p\left( S_p\geq b,X_p=a\right)- \p\left( S_p\geq b+1,X_p=a\right)\\
=\p\left( X_p=2b-a\right)- \p\left( X_p=2b+2-a\right)
=\left( \frac{1}{2}\right)^p\left[C_p^\frac{p+2b-a}{2}-C_p^{\frac{p+2b-a}{2}+1} \right]\\
=\left( \frac{1}{2}\right)^pC_p^\frac{p+2b-a}{2}\left[1-\frac{p-2b+a}{p+2b-a+2} \right]
=\p(X_p=2b-a)\frac{4b-2a+2}{p+2b-a+2}.
\end{multline*}
As $\p(X_p=2b-a)=\p(S_p=2b-a)$ (see for instance \cite{F:1} p.75) and using the remark \ref{newrem}: 
$$\frac{\p(S_p=b,X_p=a)}{\p(S_p=0)}=\frac{\p(S_p=2b-a)}{\p(S_p=0)}\frac{4b-2a+2}{p+2b-a+2} \underset{p\rightarrow\infty}{\rightarrow}0.$$
\end{proof}

\begin{lem}\label{lem6.2}
Let $y\in\mathbb N$ and $x\in\mathbb Z$. Then :
\begin{equation*}
\frac{\e\left[\mathds{1}_{y\vee(x+S_p)\neq 0}\varphi\left(y\vee(x+S_p)\right)\frac{(x+X_p)^+}{y\vee(x+S_p)}\right]}{\p(S_p=0)}
\end{equation*} 
is bounded above by $\mathds{1}_{y>0}\left\lbrace \varphi(y)\sum_{k=x}^{y}k^+  \right\rbrace+\sum_{k=y+1}^\infty k\varphi(k)$ and tends to 0 when $p\rightarrow\infty$.
\end{lem}
\begin{proof}
To simplify, we consider two cases : $\left\lbrace y>0\right\rbrace$ and $\left\lbrace y=0\right\rbrace$. For typographical simplicity we denote respectively by $\mathcal B^+$  and $\mathcal B^0$ the first and the second cases. In the first case :
\begin{eqnarray*}
\mathcal B^+:=\frac{\e\left[\varphi\left(y\vee(x+S_p)\right)\frac{(x+X_p)^+}{y\vee(x+S_p)}\right]}{\p(S_p=0)}
=\sum_{\underset{-x<\ell\leq k}{k\geq0} }\frac{\p(S_p=k,X_p=\ell)}{{\p(S_p=0)}}\varphi\left(y\vee(x+k)\right)\frac{(x+\ell)^+}{y\vee(x+k)}.
\end{eqnarray*}
According to lemma \ref{lem6.1}, we have :
\begin{eqnarray*}
\mathcal B^+&\leq&\sum_{k\geq0, -x<\ell\leq k}\varphi\left(y\vee(x+k)\right)\frac{(x+\ell)^+}{y\vee(x+k)}
\leq\sum_{k\geq0, -x<\ell\leq k}\varphi\left(y\vee(x+k)\right)\frac{(x+k)^+}{y\vee(x+k)}\\
&\leq&\sum_{k\geq0}\varphi\left(y\vee(x+k)\right)\frac{(x+k)^+(x+k)}{y\vee(x+k)}.
\end{eqnarray*}
Let us remark that we just consider cases where  $k>-x$ which implies $0< x+k\leq y\vee(x+k)$. Then :
\begin{eqnarray*}
\mathcal B^+&\leq&\sum_{k\geq 0}\varphi\left(y\vee(x+k)\right){(x+k)^+}=\varphi(y)\sum_{k=x}^{y}k^++\sum_{k=y+1}^\infty k\varphi(k).
\end{eqnarray*}
In the second case, $x\leq 0$ and $\left\lbrace y\vee(x+S_p)\neq0\right\rbrace=\left\lbrace S_p>-x\right\rbrace$. Then :
\begin{eqnarray*}
\mathcal B^0\leq \sum_{k> -x,-x<\ell\leq k} \varphi(x+k)\frac{x+\ell}{x+k}\leq \sum_{k>x}\varphi(x+k)(x+k)=\sum_{k\geq 1}k\varphi(k).
\end{eqnarray*}
We can easily conclude using lemma \ref{lem6.1} and the Lebesgue Theorem. 
\end{proof}
For $0\leq n\leq p$, one can write $S_p=S_n\vee (X_n+\tilde S_{p-n})$ where $\tilde S$ is the unilateral maximum of the standard random walk $(X_{n+k}-X_n)_{k\geq 0}$ which is independent from $\mathcal F_n$. Hence :
$$\e\left[\mathds{1}_{S_p\neq0}\varphi(S_p)\left.\frac{X_p^+}{S_p}\right| \mathcal F_n\right]=\tilde \e\left[\mathds{1}_{S_n\vee(X_n+\tilde S_{p-n})\neq0}\varphi(\tilde S_{p-n}+X_{n})\frac{(X_n+\tilde X_{p-n})^+}{(X_n+\tilde S_{p-n})\vee S_n}\right] ,$$  
where $\tilde \e$ only integrates over $\tilde S_{p-n}$ and $\tilde X_{p-n}$, $S_n$ and $X_n$ being kept fixed. Then, for $\Lambda_n\in\mathcal F_n$ :
$$\frac{\e\left[\mathds{1}_{\left\lbrace\Lambda_n,S_p\neq0\right\rbrace}\varphi(S_p)\frac{X_p^+}{S_p}\right]}{\p(S_{p-n}=0)}=\frac{\e\left[\mathds{1}_{\Lambda_n}\tilde \e\left[\mathds{1}_{S_n\vee(X_n+\tilde S_{p-n})\neq0}\varphi(\tilde S_{p-n}+X_{n})\frac{(X_n+\tilde X_{p-n})^+}{(X_n+\tilde S_{p-n})\vee S_n}\right]\right] }{\p(S_{p-n}=0)}. $$
Lemma $\ref{lem6.2}$ says that the ratio in the right hand side tends to 0 when $p$ tends to infinity and is dominated by $\varphi(S_n)\sum_{k=X_n}^{S_n} k^++\sum_{k=S_n+1}^{+\infty}k\varphi(k)$, which is integrable.\\
About the  quantity (3), we have to remark that:
\begin{equation*}
\e\left[ \mathds{1}_{y\vee(S_p+x)\neq0}(X_p+x)^+\sum_{k\geq {y\vee(S_p+x)}}\frac{\varphi(k)}{k(k+1)}\right]\leq \e\left[ \mathds{1}_{y\vee(S_p+x)\neq0}\frac{(X_p+x)^+}{y\vee (S_p+x)}\sum_{k\geq y\vee (S_p+x)}\frac{\varphi(k)}{k}\right],
\end{equation*}
and we apply the same reasoning as the one for the quantity (2) with the function $h(x)=\sum_{k\geq x}\frac{\varphi(k)}{k}$ and $x>0$ instead of $\varphi$. We just have to check that ${\sum_{x>0}x h(x)<\infty}$. Easily :
\begin{eqnarray*}
\sum_{\tilde x>0}  x h(x)=\sum  x \sum_{k\geq x}\frac{\varphi(k)}{k}
\leq \sum_{x>0}  \sum_{k\geq x}\varphi(k)
\leq\sum_{k\geq0}\sum_{k\geq x}\varphi(k)
\leq\sum_{k\geq0}k\varphi(k)<\infty.
\end{eqnarray*}
With the previous notations, we have :
\begin{equation*}
\frac{\e\left[\mathds{1}_{\Lambda_n,S_p\neq0}X_p^+\sum_{k\geq S_p}\frac{\varphi(k)}{k(k+1)}\right]}{\p(S_{p-n}=0)} \leq\frac{\e\left[\mathds{1}_{\Lambda_n}\tilde \e\left[\mathds{1}_{S_n\vee(X_n+\tilde S_{p-n})\neq0}h(\tilde S_{p-n}+X_{n})\frac{(X_n+\tilde X_{p-n})^+}{(X_n+\tilde S_{p-n})\vee S_n}\right]\right]}{\p(S_{p-n}=0)},
\end{equation*}
and we can easily conclude that the ratio in the right hand side tends to 0 when $p$ tends to infinity and is dominated by $h(S_n)\sum_{k=X_n}^{S_n} k^++\sum_{k=S_n+1}^{+\infty}kh(k)$, which is integrable.\\
To conclude the proof of the Theorem, always with the same notations we have :
$$\frac{\e\left[\mathds{1}_{\Lambda_n}f(S_p,X_p)    \right]}{\p(S_{p-n}=0)}
=\frac{\e\left[\mathds{1}_{\Lambda_n}f((S_n\vee(X_n+\tilde S_{p-n}),\tilde X_{p-n}+X_n) \right]}{\p(S_{p-n}=0)},$$
and when $p$ goes to infinity,  the ratio in the right hand side tends to $M_n^\varphi$ and is dominated by $M_n^\varphi+(\varphi(S_n)+h(S_n))\sum_{k=X_n}^{S_n} k^++\sum_{k=S_n+1}^{+\infty}k(h(k)+\varphi(k))$ which is integrable.
\setcounter{equation}{0} 
\section{Penalisation by a function $\mathbf{ S^*_{g_p }}$}
1) We start with the first point of Theorem \ref{thm8.1}. In order to prove this, we need :
\begin{lem}\label{lem7.1}
Let $\alpha>0$ and $a\in\left[-\alpha,\alpha\right]$. Then:
$$\frac{\p_a\left( S^*_{g_p}=\alpha,T_0<p\right) }{\p(S_p=0)}$$
is bounded above by 2 and tends to 1 when $p\rightarrow\infty$.
\end{lem}
To obtain this, we use a Tauberian Theorem :
\begin{thm}[Cf. \cite{F:2} p. 447]\label{thm4.2}
Given $q_n\geq0$, suppose that the series
\begin{equation*}
S(s)=\sum_{n=0}^{\infty}q_n s^n
\end{equation*}
converges for $0\leq s<1$. If $0< p<\infty$ and if the sequence 
$\left\lbrace q_n\right\rbrace $ is monotone, then the two relations:
\begin{equation*}\label{f4.1}
S(s)\underset{s\rightarrow1^-}{\sim}\frac{1}{(1-s)^p}C 
\end{equation*} 
and
\begin{equation*}
q_n\underset{n\rightarrow\infty}{\sim}\frac{1}{\Gamma(p)}n^{p-1}C,
\end{equation*}
where $0<C<\infty$, are equivalent.
\end{thm}
and the following :
\begin{lem}\label{lem7.2}
For $a<0<b$ and $\lambda\in\mathbb R $ : 
\begin{equation*}
\e\left[ \left( \cosh\lambda\right)^{-T_a\wedge T_b} \right]=\frac{\cosh\lambda\left(\frac{a+b}{2} \right)}{\cosh\lambda\left(\frac{a-b}{2} \right)}
\end{equation*}
\end{lem}
\begin{proof}
Let's recall that $X_n=\sum_{k=1}^{n}Y_k$ and define the process $\left(W_n,n\geq0 \right) $ by 
$$W_n:=\frac{\cosh\lambda\left(X_n+\beta\right)}{\left( \cosh\lambda\right)^n},$$
where $\beta\in\mathbb R$. Let's prove that $\left(W_n,n\geq0 \right) $ is a $\mathcal F_n$-martingale:
\begin{multline*}
\e\left[\cosh\lambda\left(X_{n+1}+\beta\right)\mid\mathcal F_n\right]
 =\e\left[\cosh\lambda\left(X_n+Y_{n+1}+\beta\right)\mid\mathcal F_n\right]\\
=\cosh\lambda\left(X_n+\beta\right)\e\left[ \cosh\lambda{Y_{n+1}}\right] +\sinh\lambda\left(X_n+\beta\right)\e\left[ \sin\lambda{Y_{n+1}}\right]
=\cosh\lambda\left(X_n+\beta\right)\cosh\lambda
\end{multline*}
Clearly $W$ is a martingale. Taking $\beta=-\frac{a+b}{2}$ and
using the Doob's Theorem with $T_a\wedge T_b$ :
\begin{equation}
\e\left[ W_{T_a\wedge T_b}\right]= \e\left[ W_{0}\right]=\cosh\lambda\left(\frac{a+b}{2} \right)\label{f7.3}.
\end{equation}
On the other hand, using Markov property :
\begin{eqnarray}
\e\left[ W_{T_a\wedge T_b}\right]&=&\e\left[ W_{T_a}\mathds{1}_{\left\lbrace T_a<T_b\right\rbrace} +W_{T_b}\mathds{1}_{\left\lbrace T_b<T_a\right\rbrace }\right]\nonumber\\
&=&\e\left[ \frac{\cosh\lambda\left(\frac{a-b }{2}\right)}{\left( \cosh\lambda\right)^{T_a}}\mathds{1}_{\left\lbrace T_a<T_b\right\rbrace}+\frac{\cosh\lambda\left(\frac{b-a }{2}\right)}{\left( \cosh\lambda\right)^{T_b}}\mathds{1}_{\left\lbrace T_b<T_a\right\rbrace}\right]\nonumber\\
&=&\cosh\lambda\left(\frac{a-b }{2}\right)\e\left[ \left( \cosh\lambda\right)^{-T_a\wedge T_b}\right]\label{f7.4}.
\end{eqnarray}
The formulas (\ref{f7.3}) and (\ref{f7.4}) permit to conclude.
\end{proof}
Now, we are able to prove lemma {\ref{lem7.1}} :\\
as $\p_a\left( S^*_{g_p}=\alpha,T_0<p\right)\leq \p_a\left( S_{g_p}=\alpha,T_0<p\right)$, with lemma \ref{newlem}, the first point is trivial.\\
Let $\delta_\beta$ a geometric r.v. with parameter  $0<\beta<1$ such that $\delta_\beta$ is independent of  $X$. Then:
\begin{eqnarray}
\p_a\left( S^*_{g_{\delta_\beta}}\leq\alpha\right)
=\sum_{k=1}^{\infty}\p_a\left( S^*_{g_{k}}\leq\alpha\right) \p\left( \delta_\beta=k\right)
=\sum_{k=1}^{\infty}\p_a\left( S^*_{g_{k}}\leq\alpha\right)\left( 1-\beta\right)^{k-1}\beta\label{f7.5}
\end{eqnarray}
Note that $\left\lbrace S^*_{g_p}\leq\alpha \right\rbrace=\left\lbrace g_p\leq T^*_\alpha \right\rbrace=\left\lbrace p\leq d_{T^*_\alpha}\right\rbrace=\left\lbrace p\leq T^*_\alpha+T_0.\theta_{T^*_\alpha}\right\rbrace $. Hence :
\begin{multline*}
\p_a\left( S^*_{g_{\delta_\beta}}\leq\alpha\right)=\p_a\left( \delta_\beta\leq d_{T^*_\alpha}\right)
=1-\p_a\left( \delta_\beta>d_{T^*_\alpha}\right)
=1-\e_a\left[\e_a\left[ \mathds{1}_{\delta_\beta>d_{T^*_\alpha}}\mid\mathcal F_{T^*_\alpha}\right]  \right]\\
=1-\e_a\left[ \left(1-\beta\right)^{d_{T^*_\alpha}} \right]
=1-\e_a\left[ \left(1-\beta\right)^{T^*_\alpha} \left(1-\beta\right)^{T_0.\theta_{T^*_\alpha}}\right]
=1-\e_a\left[ \left(1-\beta\right)^{T^*_\alpha}  \right]\e\left[ \left(1-\beta\right)^{ T_\alpha}\right].
\end{multline*}
We have already seen (cf. \cite{D:1} p.353  and \cite{ALR:1}):
\begin{equation*}
\e\left[ \left(1-\beta\right)^{ T_\alpha}\right]=\left( \frac{1+\sqrt{2\beta-\beta^2}}{1-\beta}\right)^{-\alpha}.
\end{equation*}
The symmetry of the quantity $\e_a\left[ \left(1-\beta\right)^{T^*_\alpha}  \right] $ permits us to assume that $a\geq0$, without a loss of generality. Then, using the Markov property and lemma \ref{lem7.2} with $(\cosh\lambda)^{-1}=1-\beta$ :
\begin{equation*}
\e_a\left[ \left(1-\beta\right)^{T^*_\alpha}  \right]=\e\left[ \left(1-\beta\right)^{T_{\left\lbrace -\alpha-a\right\rbrace}\wedge T_{\left\lbrace\alpha-a\right\rbrace }}  \right]=\frac{\cosh a\lambda}{\cosh \alpha\lambda}.
\end{equation*}
When $\beta$ goes to $0$:
\begin{equation}\label{f7.6}
\p_a\left( S^*_{g_{\delta_\beta}}\leq\alpha\right)=1-\left( \frac{1+\sqrt{2\beta-\beta^2}}{1-\beta}\right)^{-\alpha}\frac{\cosh \left[  \mathrm{argch}\left(\frac{1}{1-\beta} \right)a\right] }{\cosh\left[   \mathrm{argch}\left(\frac{1}{1-\beta}\right)\alpha\right]  }\underset{\beta\rightarrow0}{\sim}\alpha\sqrt{2 \beta}.
\end{equation}
According to the formulas (\ref{f7.5}) and (\ref{f7.6}):
\begin{equation*}
\sum_{k=1}^{\infty}\p_a\left( S^*_{g_{k}}\leq\alpha\right)\left( 1-\beta\right)^{k}\underset{\beta\rightarrow0}{\sim}\alpha\sqrt{\frac{2}{ \beta}}\left( 1-\beta\right).
\end{equation*}
In order to apply Theorem \ref{thm4.2}, put  $\beta=1-\omega$.This gives 
\begin{eqnarray*}
\sum_{k=1}^{\infty}\p_a\left( S^*_{g_{k}}\leq\alpha\right)\omega^{k}\underset{\omega\rightarrow1-}{\sim}&\alpha\omega\sqrt{\frac{2}{ \left(1-\omega \right) }}
\underset{\omega\rightarrow1-}{\sim}&\alpha\sqrt{\frac{2}{ \left(1-\omega \right) }}
\end{eqnarray*}
and this Tauberien Theorem with $p=\frac{1}{2}$ et $C= \alpha\sqrt{2}$ permits ut to obtain :
\begin{eqnarray*}
\p_a\left( S^*_{g_p}\leq\alpha\right) \underset{p\rightarrow\infty}{\sim}\frac{1}{\Gamma\left( \frac{1}{2}\right) }p^{\frac{1}{2}-1}C=\left( \frac{2}{\pi p}\right)^{\frac{1}{2}}\alpha ,
\end{eqnarray*}
and the proof can be easily finished, knowing the behaviour of $\p(S_p=0)$ when $p$ goes to $\infty$.
Thanks to this lemma, we have the following result:
\begin{lem}\label{lem7.3}
Let $x\geq 0$ and $a\in\left[-x,x\right]$. Then:
$$\frac{\e_a\left[ \varphi(x\vee S^*_{g_p})\mathds{1}_{T_0<p}\right] }{\p(S_p=0)}$$
is bounded above by $2(\varphi(x)(x-\vert a\vert )+\phi(x))$ and tends to $\varphi(x)(x-\vert a\vert )+\phi(x)$ when $p$ goes to infinity. 
\end{lem}
\begin{proof}
The proof is nearly the same as the one of lemma \ref{newlem2}
\end{proof}
With the same notations and arguments as inf Theorem \ref{thm5.1} :
\begin{equation*}
\e\left[\varphi\left( S^*_{g_p}\right)\mid\mathcal F_n  \right]
=\varphi\left( S^*_{g_n}\right)\tilde \p\left( \tilde S_{p-n}<\vert X_n\vert \right)+\e\left[\varphi\left( S^*_n\vee\tilde S^*_{\tilde g_{p-n}}\right)\mathds{1}_{\tilde T_0\leq p-n}  \right]=(1)+(2).
\end{equation*}
The end of the proof is based on the proof of Theorem \ref{thm5.1} and use lemma \ref{newlem2}. The remaining details are left to the reader. \\
Now, in order to prove that $\left(M^*_n,n\geq0 \right)$ is a martingale,  we show that conditioned by $\mathcal F_n$,   $M_{n+1}^*-M_n^*$ is zero. The case $n=0$ being trivial, we just study $n>0$.\\
First, observe that on $\left\lbrace X_n=0\right\rbrace$, $g_{n+1}=g_n=n$, $S_{n+1}^*=S_n^*$ and $\vert X_{n+1}\vert=1$. Consequently, on this event, $M^*_{n+1}-M^*_{n}=0$.\\
In the following,  we suppose that $X_n\neq 0$, and we denote $A_n:=\varphi(S_{g_{n+1}}^*)\vert X_{n+1}\vert-\varphi(S_{g_n}^*)\vert X_n\vert$ and $B_n:=\varphi(S_{n+1}^*)(S_{n+1}^*-\vert X_{n+1}\vert)-\varphi(S_n^*)(S_n^*-\vert X_n\vert)+\phi(S_{n+1}^*)-\phi(S_n^*)$. We now treat separately these two quantities :
\begin{itemize}
\item If $\left\lbrace \vert X_n\vert \geq 2\right\rbrace$,  $g_{n+1}^*=g_n^*$ and in this way $A_n=\varphi(S_{g_n}^*)(\vert X_{n+1}\vert-\vert X_n\vert)$. Conditioning on  $\mathcal F_n$, this quantity equals zero, the function  $x\rightarrow\vert x\vert$ being harmonic for the symmetric random walk except in 0.\\
If $\vert X_n\vert=1$, $A_n=2\varphi(S_{g_{n^*}})\mathds{1}_{X_{n+1}\neq 0}-\varphi(S_{g_{n^*}})$ conditional on $\mathcal F_n$ is obviously zero.\\
\item If  $\left\lbrace \vert X_n\vert\leq S_n^*-1\right\rbrace$, then $S_{n+1}^*=S_{n}^*$.
In this case,  $B_n=\varphi(S_n^*)(\vert X_{n+1}\vert-\vert X_n\vert)$ and we conclude with the harmonicity of $x\rightarrow\vert x\vert $.\\
Finally, on $\left\lbrace S_n^*=\vert X_n\vert\right\rbrace $,  $B_n=\varphi(S_{n}^*)(\mathds{1}_{S_n^*=S_{n+1}^*}-\mathds{1}_{S_n^*+1=S_{n+1}^*})$ and conditionned on $\mathcal F_n$, it is clear that this quantity equals zero.
\end{itemize}
Consequently $M^*$ is a martingale satisfying :
$$\left|M_n^*-M_0^* \right|\leq 3n ,$$
and as $M_0^*=1$, one has $\e[M_n^*]=1$. Observe that the positivity of $M^*$ is obvious from the definitions of $\varphi$ and $\phi$. \\
 2) Now, we prove point 2 of Theorem \ref{thm8.1}. 
For $a>0$:
\begin{equation*}
Q^*\left( g_p>a\right) 
=\e^{Q^*}\left[\mathds{1}_{p>d_a} \right]
=\e\left[\mathds{1}_{p>d_a}M^*_{d_a} \right]
=\e\left[\mathds{1}_{p>d_a}\left\lbrace \varphi\left( S_{d_a}^*\right) S_{d_a}^* +\phi\left( S_{d_a}^*\right)\right\rbrace \right].
\end{equation*}
As $a$ is fixed, the sequence of positive random variables $\left(\mathds{1}_{p>d_a}\left\lbrace \varphi\left( S_{d_a}^*\right) S_{d_a}^* +\phi\left( S_{d_a}^*\right)\right\rbrace\right)_{p\geq0}$ is increasing and tends to $ \varphi\left( S_{d_a}^*\right) S_{d_a}^* +\phi\left( S_{d_a}^*\right)$, and  the sequence of events $\left\lbrace g_p>a\right\rbrace$ is increasing and tends to $\left\lbrace g>a\right\rbrace$ when $p$ tends to infinity . Hence, according to Lebesgue Theorem,  when  $p$ goes to $+\infty$:
\begin{equation*}
Q^*\left( g>a\right)=\e\left[ \varphi\left( S_{d_a}^*\right) S_{d_a}^* +\phi\left( S_{d_a}^*\right) \right].
\end{equation*}
As $\phi(S_{d_a}^*)\leq 1$,  Lebesgue Theorem implies that $\e\left[\phi\left( S_{d_a}^*\right) \right]\underset{a\rightarrow\infty}{\rightarrow}0$.\\
It remains to prove $\e\left[ \varphi\left( S_{d_a}^*\right) S_{d_a}^*\right] \underset{a\rightarrow\infty}{\rightarrow}0$.
\begin{lem}\label{lem7.4}
Let $\psi:\mathbb N\rightarrow\mathbb R^+$ such that $\sum_{k\geq 0}\varphi(k)<+\infty$. For $a>0$ :
\begin{equation*}
\e\left[ \psi\left( S^*_{d_a}\right) \mid\mathcal F_{a}\right]=\mathds{1}_{X_a=0}\psi\left( S_a^*\right)+\mathds{1}_{X_a\neq0}\left\lbrace 
\psi\left(S^*_a \right) \left(1-\frac{\vert X_a\vert}{S^*_a}  \right)+\vert X_a\vert\sum_{k=S_a}^\infty \frac{\psi(k)}{k(k+1)}\right\rbrace .
\end{equation*}
\end{lem}
\begin{proof}
Let $\tilde S^*_{\tilde T_0}$ be the bilateral maximum of a walk issued from $X_a$ until the hitting time of the level $0$ and which is independent of $\mathcal F_a$.
\begin{eqnarray*}
\e\left[ \psi\left( S^*_{d_a}\right) \mid\mathcal F_{a}\right]&=&\e\left[ \psi\left( S^*_{a}\vee\tilde S^*_{\tilde T_0}\right) \mid\mathcal F_{a}\right]=\tilde \e_{\vert X_a\vert}\left[ \psi\left( S^*_{a}\vee\tilde S^*_{\tilde T_0}\right)\right]\\
&=&\mathds{1}_{X_a=0}\psi\left( S^*_{a}\right)
+\mathds{1}_{X_a\neq0}\sum_{k\geq\vert X_a\vert}\p_{X_a}\left( \tilde S^*_{\tilde T_0}=k\right)\psi\left(S^*_{a}\vee k \right) .
\end{eqnarray*}
On $\left\lbrace X_a>0\right\rbrace$, as the sign of $\tilde X$ does not change between 0 and $\tilde T_0$,  $\left\lbrace  \tilde S^*_{\tilde T_0}=k\right\rbrace=\left\lbrace  \tilde S_{\tilde T_0}=k\right\rbrace$.
Moreover, thanks to the symmetry of $\tilde X$, on $\left\lbrace X_a<0\right\rbrace$, $\left\lbrace \tilde S^*_{\tilde T_0}=k\right\rbrace= \left\lbrace\tilde S_{\tilde T_0}=k\right\rbrace$.  So, according to lemma \ref{lem5.5}:
\begin{equation*}
\mathds{1}_{X_a\neq0}\p_{X_a}\left( \tilde S^*_{\tilde T_0}=k\right)=\mathds{1}_{X_a\neq0}\p_{X_a}\left( \tilde S_{\tilde T_0}=k\right)=\frac{\vert X_a\vert}{k(k+1)}.
\end{equation*}
Consequently on $\left\lbrace X_a\neq0\right\rbrace $ :
\begin{eqnarray*}
\sum_{k\geq\vert X_a\vert}\p_{X_a}\left( \tilde S^*_{\tilde T_0}=k\right)\psi\left(S^*_{a}\vee k \right) &=&\sum_{k\geq\vert X_a\vert}\psi\left(S^*_{a}\vee k \right) \frac{\vert X_a\vert}{k(k+1)}\\
&=&\sum_{k=\vert X_a\vert}^{S^*_{a}-1}\psi\left(S^*_{a}\right) \frac{\vert X_a\vert}{k(k+1)}+\sum_{k\geq S^*_{a}}\psi\left( k \right) \frac{\vert X_a\vert}{k(k+1)}\\
&=&\psi\left(S^*_{a}\right)\vert X_a\vert\left(\frac{1}{\vert X_a\vert}-\frac{1}{S^*_a} \right)+\sum_{k\geq S^*_{a}}\psi\left( k \right) \frac{\vert X_a\vert}{k(k+1)}.
\end{eqnarray*}
\end{proof}
Applying this lemma with $\psi(x)=x\varphi(x)$ :
\begin{eqnarray*}
\e\left[ \varphi\left( S_{d_a}^*\right) S_{d_a}^*\right]
&=&\e\left[ \mathds{1}_{X_a=0}\varphi\left( S^*_a\right)S^*_a+ \mathds{1}_{X_a\neq0}\left\lbrace\varphi\left( S^*_a\right)\left(S^*_a-\vert X_a\vert \right)+\vert X_a\vert\sum_{k\geq S^*_a}\frac{\varphi(k)}{k+1}\right\rbrace\right]\\
&\leq&\e\left[ \varphi\left( S^*_a\right)S^*_a\right]+\e\left[ \frac{\vert X_a\vert}{S^*_a+1}\sum_{k\geq S^*_a}\varphi(k)\right]
\leq\e\left[ \varphi\left( S^*_a\right)S^*_a\right]+\e\left[ \sum_{k\geq S^*_a}\varphi(k)\right]\\
&\leq&\e\left[ \varphi\left( S^*_a\right)S^*_a\right]+\e\left[ \phi\left( S^*_a\right) \right].
\end{eqnarray*}
As $\phi(S^*_a)\leq 1$ and $\phi(S_a^*)$ tends to 0 a.s. when $a$ tends to infinity, Lebesgue Theorem implies that $\e[\phi(S_a^*)]\underset{a\rightarrow+\infty}{\rightarrow}0$ .
On the other hand:
\begin{eqnarray*}
\e\left[ \varphi\left( S^*_a\right)S^*_a\right]=\sum_{k\geq0}\varphi(k)k\p\left(S^*_a=k \right)\leq\sum_{k\geq0}\varphi(k)k\p\left(S_a=k \right)
\leq \e\left[ \varphi\left( S_a\right)S_a\right],
\end{eqnarray*}
and we have already proved that $\e[\varphi(S_a)S_a]$ tends to 0 when $a$ tends to infinity (cf. point 2 Theorem \ref{thm5.1}).
As a result $g$ is $Q$-a.s. finite and :
\begin{equation*}
Q^*\left( g=\infty\right)=\lim_{a\rightarrow\infty}Q^*\left( g>a\right)=0.
\end{equation*}
\hfill\qed\\
3) We now prove the third and last point of the Theorem.
\begin{lem}\label{lem7.5}
For all $a>1$:
$$\p\left(S^*_{\tau_a}=k \right)= \left\lbrace
\begin{array}{cl}
0&\mbox{, if $k=0$}\\
\left( 1-\frac{1}{k+1}\right)^{a-1} -\left( 1-\frac{1}{k}\right)^{a-1}&\mbox{, otherwise.}
\end{array}
\right.$$
\end{lem}
\begin{proof}
Using Markov property and thanks to the symmetry of $X$ :
\begin{eqnarray*}
\p\left(\tau_2<T^*_k \right)&=&\frac{1}{2}\left[ \p_1\left(T_0<T^*_k \right)+\p_{-1}\left(T_0<T^*_k \right)\right]\\
&=&  \frac{1}{2}\left[ \p_1\left(T_0<T_k \right)+\p_{-1}\left(T_0<T_{-k} \right)\right]=\p_1\left(T_0<T_k \right).
\end{eqnarray*}
Recall that $\p_1(T_0<T_k)=1-\frac{1}{k}$.
Moreover using the strong Markov property and an obvious recurrence :
\begin{eqnarray*}
\p\left(S^*_{\tau_a}<k \right)&=&\p\left( \tau_a<T^*_k\right)=\p\left(\tau_a<T^*_k\vert \tau_2<T^*_k \right)\p\left(\tau_2<T^*_k\right)\\
&=&\p\left(\tau_{a-1}<T^*_k \right)\p\left(\tau_2<T^*_k\right)=\p\left(\tau_{2}<T^*_k \right)^{a-1}=\left(1-\frac{1}{k} \right)^{a-1}.
\end{eqnarray*}
\end{proof}

Using a similar reasoning and notations as the one of $\varphi(S_{g_p})$  :
\begin{equation*}
M^*_{\tau_a}-M^*_{\tau_{a+1}}=\mathds{1}_{\tilde S^*_{\tilde \tau_2}>S^*_{\tau_a}}(\varphi(S^*_{\tau_a})S^*_{\tau_a}-\varphi(\tilde S^*_{\tilde \tau_2})\tilde S^*_{\tilde \tau_2}+\sum_{k=S^*_{\tau_a}}^{\tilde S^*_{\tilde \tau_2}-1}\varphi(k))
\Rightarrow 
\e\left[ M_{\tau_a}-M_{\tau_{a+1}}\vert \mathcal F_{\tau_a}\right]=\varphi(S^*_{\tau_a}).
\end{equation*}
Let  $F$ be a positive functional,  $f_1$ and $f_2$ be two functions from $\mathbb N$ to $\mathbb R^+$. :
\begin{equation*}
\mathcal G:=\e^{Q^*}\left[F\left( X_u,u\leq g\right)f_1\left(\gamma_g \right)f_2\left( S^*_g\right)    \right]=\sum_{a\geq1}\e\left[F\left( X_u,u\leq \tau_a\right)f_1\left(a \right)f_2\left( S^*_{\tau_a}\right)\varphi\left(S^*_{\tau_a} \right)    \right].
\end{equation*}
With $F\equiv1$:
\begin{eqnarray*}
\mathcal G&=&\sum_{a\geq0}\sum_{k\geq0}f_1(a)f_2(k)\varphi(k)\p\left( S^*_{\tau_a}=k\right)\\
&=&\mathds{1}_{k=0,a=1}\varphi(0)f_1(1)f_2(0)+\sum_{a>1,k>0}\varphi(k)\left[\left(1-\frac{1}{k+1} \right)^{a-1}-\left(1-\frac{1}{k} \right)^{a-1} \right]f_1(a)f_2(k).
\end{eqnarray*}
Then, the law of $\left( \gamma_g,S^*_g\right) $ is:
\begin{equation*}
Q^*\left( \gamma_g=a,S_g^*=k\right) =\mathds{1}_{k=0,a=1}\varphi(0)+\mathds{1}_{k>0,a>1}\varphi(k)\left[\left(1-\frac{1}{k+1} \right)^{a-1}-\left(1-\frac{1}{k} \right)^{a-1} \right]
\end{equation*}
We easily find the density of $S^*_g$ summing over $a$.\\
Writing  $\mathcal G$ in two different ways :
\begin{eqnarray*}
\mathcal G
&=&\sum_{a\geq1}f_1\left(a\right)f_2\left( k\right)\varphi\left( k\right) Q^*(\gamma_g=a,S^*_g=k)\e^{Q^*}\left[F\left( X_u,u\leq \tau_a\right)    \right]\\
&=& \sum_{a\geq1}\sum_{k\geq0}f_1(a)f_2(k)\p\left( S^*_{\tau_a}=k\right)\varphi(k)E\left[ F\left( X_u,u\leq\tau_{a}\right)\mid S^*_{\tau_a}=k \right] ,
\end{eqnarray*}
we conclude that :
\begin{equation*}
\e^Q\left[F\left( X_u,u\leq g\right)\mid S^*_g=k,\gamma_g=a  \right]=E\left[ F\left( X_u,u\leq\tau_{a}\right)\mid S^*_{\tau_a}=k \right].
\end{equation*}
This achieves the proof of point 3.iii. \\
3.iii) The  study of the process $\left(X_u,u\leq g \right)$ under $Q^*$ is very close to the study of $(X_n,n\geq0)$ under $Q$  in Theorem $\ref{thm5.1}$ :
\begin{lem}\label{lem7.6}
Under the law $Q^*$ and  conditional on 
the event $\mathcal H_l$, the
process $(X_{n+\tau_l},\ {1\leq n\leq T_{p,\tau_l}-\tau_l})$
 is a 3-Bessel* walk started from 1 and stopped when it first hits the level $p$.\end{lem}
\begin{proof}
We just have to see that $M_{\tau_{l+1}}^*=\varphi(S^*_{\tau_{l+1}})S^*_{\tau_{l+1}}+\phi(S^*_{\tau_{l+1}})$ and conditioning by $\mathcal F_{T_{p,\tau_l}}$, we obtain $\e\left[M_{\tau_{l+1}}^*\vert \mathcal F_{T_{p,\tau_l}} \right] =\e_p \left[  \varphi(S^*_{T_0})S^*_{T_0}+\phi(S^*_{T_0}) \right]$ is  a constant. 
\end{proof}
We can easily prove by symmetry that $(X_{n+g},n\geq 0)$ is either a 3-dimensional Bessel* walk either a reversed 3-Bessel* walk. It remains to know with what probability we have each case. 
Nevertheless we can deduce from the previous result that under $Q^*$, $S^*_\infty=\infty$. It permits us to obtain the following lemma :
\begin{lem}\label{lem7.7}
Under $Q^*$, $S^*_{g_{T_a}}$ is a uniformly distributed random variable on $\left\lbrace 0,1,\dots,a-1\right\rbrace$.
\end{lem}
\begin{proof}
According to Doob's Theorem :
\begin{eqnarray*}
&Q^*\left( S^*_p>a\right)&=Q^*\left( T^*_a<p\right)
=\e\left[\mathds{1}_{T^*_a<p}M^*_{T^*_a} \right]
=\e\left[M^*_{T^*_a} \right]=\e\left[\mathds{1}_{T^*_a<p}\varphi\left(S^*_{g_{T^*_a}} \right)a+\phi\left( a\right) \right].
\end{eqnarray*}
When $p$ tends to infinity :
$$1=Q^*\left(S^*_\infty>a \right)=\e\left[\varphi\left(S^*_{g_{T^*_a}} \right)a+\phi\left( a\right) \right] \Leftrightarrow \sum_{k=0}^{a-1}\varphi(k)=a\sum_{k=0}^{a-1}\p\left(S^*_{g_{T^*_a}}=k\right)\varphi(k).$$
The fact this equality is true for a family of function $\varphi$ (for example $\varphi_\lambda(x)=e^{-\lambda x}$) permits us to say that $\forall k\in\left\lbrace 0,1,\dots,a-1\right\rbrace$, $\p\left(S^*_{g_{T^*_a}}=k\right)=a^{-1}$
\end{proof}
Recall $\Delta^+:=\left\lbrace X_{n+g}>0,\forall n>0\right\rbrace$ (resp.  $\Delta^-:=\left\lbrace X_{n+g}<0,\forall n>0\right\rbrace$). As $g$ is $Q^*$-a.s. finite,
$Q^*\left(\Delta^+ \right)=\lim_{p\rightarrow\infty}Q^*\left(X_{T^*_p}>0 \right)$ and with the definition of $Q^*$ :
\begin{eqnarray*}
 Q^*\left(X_{T^*_p}>0 \right)=\e\left[\mathds{1}_{X_{T_p^*}>0}M_{T^*_p}\right]=\e\left[\mathds{1}_{X_{T_p^*}>0}\left\lbrace  \varphi\left(S^*_{g_{T_p^*}}\right)p+\phi(p)  \right\rbrace\right].
\end{eqnarray*}
Using the symmetry of $X$ under $\p$ :
$$2\e\left[\mathds{1}_{X_{T_p^*}>0}\varphi\left(S^*_{g_{T_p^*}}\right)\right]=\e\left[\mathds{1}_{X_{T_p^*}>0}\varphi\left(S^*_{g_{T_p^*}}\right)\right]+\e\left[\mathds{1}_{X_{T_p^*}<0}\varphi\left(S^*_{g_{T_p^*}}\right)\right]=\e\left[\varphi\left(S^*_{g_{T_p^*}}\right)\right]=\sum_{k=0}^{p-1}\frac{\varphi(k)}{p}.$$ 
Consequently $Q^*\left(X_{T^*_p}>0 \right)=
\frac{1}{2}$
and :
$$Q^*\left(\Delta^+   \right)=\lim_{p\rightarrow\infty}Q^*\left(X_{T^*_p}>0\right)=\frac{1}{2}.$$

\section{Penalisation by  $S^*_p$}
In fact, we have a better result :
\begin{thm}
\begin{enumerate}
\item Let $a,b>0$, then :
\begin{equation}\label{fminmax2}
\lim_{p\rightarrow\infty}\frac{\e\left[\mathds{1}_{\left\lbrace\Lambda_n,\,S_p<a,I_p>-b\right\rbrace}\right]}{\e\left[\mathds{1}_{\left\lbrace S_p<a,I_p>-b\right\rbrace}\right]}:=\e\left[\mathds{1}_{\left\lbrace\Lambda_n,S_n<a,I_n>-b\right\rbrace}M_n\right],
\end{equation}
where $M_n:=\left(\cos\left({\frac{\pi}{a+b}}\right)\right)^{-n}\frac{\sin\left( \frac{\pi(a-X_n)}{a+b}\right)}{\sin\left( \frac{\pi a}{a+b}\right)}$ is a positive martingale non uniformly integrable.
\item Let us define a new probability $Q$ on $\left(\Omega,\,\mathcal F_\infty\right)$ characterized by :
\begin{equation}
\forall n\in\mathbb N,\,\forall \Lambda_n\in\mathcal F_n,\, Q\left( \Lambda_n\right):=\e\left[\Lambda_nM_n\right].
\end{equation}
Then under $Q$, $\left(X_n,n\geq0\right)$ have the following transition probabilities for 
$-b+1\leq k\leq a-1$:
\begin{eqnarray*}
Q\left(X_{n+1}=k+1\vert X_n=k\right)&=&\frac{\sin\left(\frac{a-k-1}{a+b}\pi\right)}{2\cos\left(\frac{\pi}{a+b}\right)\sin\left(\frac{a-k}{a+b}\pi\right)},\\
Q\left(X_{n+1}=k-1\vert X_n=k\right)&=&\frac{\sin\left(\frac{a-k+1}{a+b}\pi\right)}{2\cos\left(\frac{\pi}{a+b}\right)\sin\left(\frac{a-k}{a+b}\pi\right)}.
\end{eqnarray*}

\end{enumerate}
\end{thm}
1) To prove the first point of Theorem we need the following lemma :
\begin{lem}\label{leminmax}
Let $a,b>0$ and  $c\in\left[-b+1,a-1 \right]$, then:
\begin{equation}\label{formulfinish}
\p\left(S_n<a,I_n>-b\right)\underset{n\rightarrow\infty}{\sim}\frac{4}{a+b}\left(\cos\left(\frac{\pi}{a+b}\right)\right)^n\sin\left(\frac{a\pi}{a+b}\right)\sum_{\underset{c\equiv n[2]}{c=-b+1}}^{a-1}\sin\left(\frac{\pi(a-c)}{a+b}\right).
\end{equation}
\end{lem}
Let us postpone the proof of this lemma and finish the proof of (\ref{fminmax2}).
As usual, let $\tilde X_k=X_{k+n}$, a random walk started from $X_n$ and independent of $\mathcal F_n$, and $\tilde S_n$ and $\tilde I_n$ respectively the supremum and infimum associated to $\tilde X$. In the following steps $\tilde \p$ is the measure associated to  $\tilde X$, $X_n$, $S_n$ and $I_n$ being kept fixed. Using the Markov property : 
\begin{eqnarray*}
\e\left[\mathds{1}_{\left\lbrace\Lambda_n,\,S_p<a,I_p>-b\right\rbrace}\right]&=&\e\left[\mathds{1}_{\left\lbrace\Lambda_n,\,S_n<a,I_n>-b\right\rbrace}\tilde \p\left( \tilde S_{p-n}<a-X_n,\tilde I_{p-n}>-b-X_n\right)\right].
\end{eqnarray*}
 Lemma \ref{leminmax} says:
\begin{multline*}
 \p\left( \tilde S_{p-n}<a-X_n,\tilde I_{p-n}>-b-X_n\right)\\
\underset{p\rightarrow\infty}{\sim}
\sum_{c=-b-X_n+1,c\equiv p-n\,\left[2\right]}^{a-X_n-1}\frac{4}{a+b}\left(\cos\left({\frac{\pi}{a+b}}\right)\right)^{p-n}\sin\left(\frac{\pi (a-X_n)}{a+b}\right)\sin\left( \frac{\pi(a-X_n-c)}{a+b}\right)\\
\underset{p\rightarrow\infty}{\sim}
\frac{4}{a+b}\left(\cos\left({\frac{\pi}{a+b}}\right)\right)^{p-n}\sin\left(\frac{\pi (a-X_n)}{a+b}\right)\sum_{c=-b+1,c\equiv p\,\left[2\right]}^{a-1}\sin\left( \frac{\pi(a-c)}{a+b}\right).
\end{multline*}
Dividing this formula by (\ref{formulfinish}), we obtain (\ref{fminmax2}).\\
\textbf{Proof of lemma \ref{leminmax} :}\\
To prove this lemma we need the following combinatory result : 
\begin{lem}\label{lemcombi}
Let $p\in\mathbb N$, $0<u<p$ :
$$\sum_{k\geq 0} C_n^{kp+u}=\frac{1}{p}\sum_{\ell=0}^{p-1}\left( 1+e^{\frac{\ell 2i\pi}{p}} \right)^ne^{-\frac{2i\pi\ell u}{u}}.  $$
\end{lem}
\begin{proof}[Proof of lemma \ref{lemcombi}:  ]
\begin{equation*}
\sum_{\ell=0}^{p-1}\left(1+e^{\frac{\ell2i\pi}{p}}\right)^ne^{-\frac{2i\pi\ell u}{p}}=\sum_{\ell=0}^{p-1}\sum_{k=0}^nC_n^ke^{\frac{2i\pi\ell k}{p}}e^{-\frac{2i\pi\ell u}{p}}
\sum_{k=0}^nC_n^k\sum_{\ell=0}^{p-1}e^{\frac{2i\pi\ell (k-u)}{p}}.
\end{equation*}
Those sums can be easily simplifed if we note that  :
$$\sum_{\ell=0}^{p-1}e^{\frac{2i\pi\ell (k-u)}{p}}=\left\lbrace
\begin{array}{cl}
p&\mbox{, if $k-u$ is a multiple of $p$}\\
\sum_{\ell=0}^{p-1}e^{\frac{2i\pi\ell (k-u)}{p}}=\frac{1-e^{2i\pi(k-u)}}{1-e^{\frac{2i\pi(k-u)}{p}}}=0&\mbox{, otherwise.}
\end{array}
\right. $$
Then :
$$\sum_{\ell=0}^{p-1}\left(1+e^{\frac{\ell2i\pi}{p}}\right)^ne^{-\frac{2i\pi\ell u}{p}}=\sum_{k=0,k\equiv u\left[p \right] }^npC_n^k
=p\sum_{k\geq0}C_n^{kp+u}.$$
\end{proof}
According to \cite{F:1} p.79:
\begin{equation*}
\p\left(S_n<a,X_n=c,I_n>-b\right)=\left(\frac{1}{2}\right)^n\sum_{k\in\mathbb Z}C_n^{\frac{n+c}{2}+k(a+b)}-C_n^{\frac{n-c}{2}+k(a+b)+a}.
\end{equation*}
Clearly $c$ and $n$ must have the same parity and denote :
$$A_n^c=\nicefrac{(n+c)}{2}+k_0(a+b),
B_n^c=\nicefrac{(n-c)}{2}+a+k_1(a+b),$$
where $k_0$ (respectively $k_1$) is the first  $k$ such as $\nicefrac{(n+c)}{2}+k(a+b)$ (respectively $\nicefrac{(n-c)}{2}+k(a+b)+a$ ) is positive.
Then :
\begin{multline*}
\p\left(S_n<a,X_n=c,I_n>-b\right)=\frac{2^{-n}}{a+b}\sum_{\ell=0}^{a+b-1}\left(1+e^{\frac{2i\pi\ell}{a+b}}\right)^{n}\left[e^{-\frac{2i\pi\ell A_n^c}{a+b}}-e^{-\frac{2i\pi\ell B_n^c}{a+b}}\right]\\
=\frac{2i}{a+b}\sum_{\ell=1}^{a+b-1}\cos^n\left({\frac{\pi\ell}{a+b}}\right)e^{\frac{i\pi\ell(n-(A_n^c+B_n^c))}{a+b}}\sin\left(\frac{\pi\ell(B_n^c-A_n^c)}{a+b}\right)\\
=-\frac{2}{a+b}\sum_{\ell=1}^{a+b-1}\cos^n\left({\frac{\pi\ell}{a+b}}\right)\sin{\frac{\pi\ell( n-(A_n^c+B_n^c))}{a+b}}\sin\left(\frac{\pi\ell(B_n^c-A_n^c)}{a+b}\right).
\end{multline*}

Let us remark that:
$$
\left\lbrace
\begin{aligned}
B_n^c-A_n^c=-a-(k_0+k_1)(a+b)\\
n-A_n^c-B_n^c=a-c+(k_1-k_0)(a+b)
\end{aligned}\right.
$$
Which implies :
\begin{multline*}
\p\left(S_n<a,X_n=c,I_n>-b\right)=\\
\frac{2}{a+b}\sum_{\ell=1}^{a+b-1}(-1)^{(k_0+k_1)\ell+(k_1-k_0)\ell}\cos^n\left({\frac{\pi\ell}{a+b}}\right)\sin{\frac{\pi\ell a}{a+b}}\sin\left(\frac{\pi\ell(a-c))}{a+b}\right)\\
=\frac{2}{a+b}\sum_{\ell=1}^{a+b-1}\cos^n\left({\frac{\pi\ell}{a+b}}\right)\sin{\frac{\pi\ell a}{a+b}}\sin\left(\frac{\pi\ell(a-c))}{a+b}\right).
\end{multline*}
Here, we have to notice that when $n\rightarrow\infty$, the leading terms are $\ell=1$ and $\ell=a+b-1$. Hence :
\begin{multline*}
\cos^n\left({\pi -\frac{\pi}{a+b}}\right)\sin\left({\pi a-\frac{\pi a}{a+b}}\right)\sin\left(\pi(a-c)-\frac{\pi(a-c))}{a+b}\right)=\\
(-1)^{n+2a-c}\cos^n\left({\frac{\pi}{a+b}}\right)\sin\left(\frac{\pi a}{a+b}\right)\sin\left( \frac{\pi(a-c)}{a+b}\right)
=\cos^n\left({\frac{\pi}{a+b}}\right)\sin\left(\frac{\pi a}{a+b}\right)\sin\left( \frac{\pi(a-c)}{a+b}\right),
\end{multline*}
$n$ and $c$ having the same parity. 
Hence:
\begin{equation*}
\p\left(S_n<a,X_n=c,I_n>-b\right)\underset{n\rightarrow\infty}{\sim}\frac{4}{a+b}\cos^n\left({\frac{\pi}{a+b}}\right)\sin\left(\frac{\pi a}{a+b}\right)\sin\left( \frac{\pi(a-c)}{a+b}\right)
\end{equation*}
and we can easily deduce :
\begin{equation*}
\p\left(S_n<a,I_n>-b\right)\underset{n\rightarrow\infty}{\sim}\frac{4}{a+b}\cos^n\left({\frac{\pi}{a+b}}\right)\sin\left(\frac{\pi a}{a+b}\right)\sum_{c\equiv n[2],c=-b+1}^{a-1}\sin\left( \frac{\pi(a-c)}{a+b}\right).
\end{equation*}
We need to prove that $M$ is a positive martingale. Positivity is obvious and for all $n\geq 0$  :
$$M_n \leq \frac{\left(\cos\left({\frac{\pi}{a+b}}\right)\right)^{-n}}{\sin\left( \frac{\pi a}{a+b}\right)}$$
\begin{eqnarray*}
\e\left[\left.\sin\left(\frac{\pi(a-X_{n+1})}{a+b}\right)\right| \mathcal F_n\right]=\e\left[\left.\sin\left(\frac{\pi(a-X_{n})}{a+b}\right)\cos\left(\frac{\pi Y_{n+1}}{a+b}\right)+\right|\mathcal F_n\right]\\+\e\left[\left.\cos\left(\frac{\pi(a-X_{n})}{a+b}\right)\sin\left(\frac{\pi Y_{n+1}}{a+b}\right)\right|\mathcal F_n\right]
=\sin\left(\frac{\pi(a-X_{n})}{a+b}\right)\cos\left(\frac{\pi }{a+b}\right)
\end{eqnarray*}
Hence $M$ is a martingale stopped when it hits the boundary of the segment $\left[-b,a\right]$.\\
2) For $-b+1\leq k\leq a-1$, using the Markov property and the defintion of $Q$ :
\begin{multline*}
Q\left(X_{n+1}=k+1\vert X_n=k \right)=\frac{Q\left(X_{n+1}=k+1,\,X_n=k \right)}{Q\left(X_n=k\right)}
=\frac{\e\left[\mathds{1}_{\left\lbrace X_{n+1}=k+1,\,X_n=k \right\rbrace}M_{n+1}\right]}{\e\left[\mathds{1}_{\left\lbrace X_n=k \right\rbrace}M_{n}\right]}\\
=\frac{\e\left[\mathds{1}_{\left\lbrace X_{n+1}=k+1,\,X_n=k , S_{n+1}<a,I_{n+1}>-b\right\rbrace}\left(\cos\left({\frac{\pi}{a+b}}\right)\right)^{-n-1}\frac{\sin\left( \frac{\pi(a-k-1)}{a+b}\right)}{\sin\left( \frac{\pi a}{a+b}\right)}\right]}{\e\left[\mathds{1}_{\left\lbrace X_n=k ,\Lambda_n,S_n<a,I_n>-b\right\rbrace}\left(\cos\left({\frac{\pi}{a+b}}\right)\right)^{-n}\frac{\sin\left( \frac{\pi(a-k)}{a+b}\right)}{\sin\left( \frac{\pi a}{a+b}\right)}\right]}\\
=\frac{\left(\cos\left({\frac{\pi}{a+b}}\right)\right)^{-1}\sin\left( \frac{\pi(a-k-1)}{a+b}\right)}{\sin\left( \frac{\pi(a-k)}{a+b}\right)}\p\left({ X_{n+1}=k+1\vert X_n=k,\, S_n<a,\,I_n>-b }\right)\\
=\frac{\left(\cos\left({\frac{\pi}{a+b}}\right)\right)^{-1}\sin\left( \frac{\pi(a-k-1)}{a+b}\right)}{2\sin\left( \frac{\pi(a-k)}{a+b}\right)}
\end{multline*}

\begin{rem}
We can easily complete this study by a penalisation functional $G_p=\mathds{1}_{S^*_{d_p}<a}$. We just have to see that this penalisation is the same as $\mathds{1}_{S^*_p<a}$.
\end{rem}

\nocite{RVY:1}\nocite{F:1}\nocite{LG:1}\nocite{ALR:1}\nocite{F:2}\nocite{D:1}
\providecommand{\bysame}{\leavevmode\hbox to3em{\hrulefill}\thinspace}
\providecommand{\MR}{\relax\ifhmode\unskip\space\fi MR }
\providecommand{\MRhref}[2]{%
  \href{http://www.ams.org/mathscinet-getitem?mr=#1}{#2}
}
\providecommand{\href}[2]{#2}


\begin{thebibliography}{RVY06}

\bibitem[ALR04]{ALR:1}
{C}. Ackermann, {G}. {L}orang, and {B}. {R}oynette, \emph{Independance of time
  and position for a random walk}, Revista {M}atematica {I}beroamericana
  \textbf{20} (2004), no.~3, pp. 915--917.

\bibitem[Deb09]{D:1}
{P}ierre Debs, \emph{Penalisation of the standard random walk by a function of
  the one-sided maximum, of the local time, or of the duration of the
  excursions}, S\'eminaire de probabilit\'es XLII (2009).

\bibitem[Fel50]{F:1}
Feller, \emph{An {I}ntroduction to {P}robability {T}heory and {I}ts
  {A}pplications}, vol.~1, 1950.

\bibitem[Fel71]{F:2}
\bysame, \emph{An {I}ntroduction to {P}robability {T}heory and {I}ts
  {A}pplications}, vol.~2, 1966-1971.

\bibitem[LeG85]{LG:1}
{J}.{F}. LeGall, \emph{Une {A}pproche {E}l\'ementaire des {T}h\'eor\`emes de
  {D}\'ecomposition de {W}illiams}, Lecture {N}otes in {M}athematics,
  {S}\'eminaire de {P}robabilit\'es XX, 1984-1985, pp.~447--464.

\bibitem[RVY]{RVY:7}
{B}. Roynette, {P}. {V}allois, and {M}. {Y}or, \emph{Brownian penalisations
  related to excursion lengths}, to be published in {A}nnales de l'Institut
  {H}enri {P}oincar\'e.

\bibitem[RVY06]{RVY:1}
\bysame, \emph{Limiting {L}aws associated with {B}rownian {M}otion perturbed by
  {I}ts {M}aximum, {M}inimum and {L}ocal {T}ime}, Studia sci. {H}ungarica
  {M}athematica \textbf{43} (2006), no.~3.

\end{thebibliography}
\end{document}